\newcount\ite\ite=1\def\0{\global\ite=1\1}
\def\1{\item{\rm(\romannumeral\the\ite)}\advance\ite1\quad}

\documentclass[12pt]{amsart}

\pdfoutput=1

\usepackage{epic, eepic, amsfonts, latexsym, amssymb, graphicx,
multicol, mathrsfs, color, amscd, verbatim, paralist,
xspace, url, euscript, stmaryrd,  amsmath, enumitem,
bbold, multirow, tikz}
\usepackage[all]{xy}

\usepackage[colorlinks, linkcolor=blue, citecolor=magenta, urlcolor=cyan]{hyperref}

\font\teneufm=eufm10 scaled \magstep1
\font\seveneufm=eufm7 scaled \magstep1
\font\fiveeufm=eufm5  scaled \magstep1
\newfam\eufmfam

\textfont\eufmfam=\teneufm
\scriptfont\eufmfam=\seveneufm
\scriptscriptfont\eufmfam=\fiveeufm

\newfam\msbfam
\font\tenmsb=msbm10 scaled \magstep1  \textfont\msbfam=\tenmsb
\font\sevenmsb=msbm7 scaled \magstep1 \scriptfont\msbfam=\sevenmsb
\font\fivemsb=msbm5 scaled \magstep1  \scriptscriptfont\msbfam=\fivemsb

\def\dd#1{\raise1.5pt\hbox{$\,\partial\!$}/\raise-2.5pt\hbox{$\!\partial#1\,$}}

\def\tilde{\widetilde}
\def\hat{\widehat}
\def\5#1{{\mathcal #1}}

\def\RR{{\mathbb R}}
\def\CC{{\mathbb C}}
\def\BB{{\mathbb B}}

\def\NN{{\mathbb N}}
\def\ZZ{{\mathbb Z}}
\def\TT{{\mathbb T}}
\def\PP{{\mathbb P}}
\def\HH{{\mathbb H}}

\def\ra{\rightarrow}

\def\hat{\widehat}

\def\GL{\mathop{\rm GL}\nolimits}

\def\Im{\mathop{\rm Im}\nolimits}
\def\Re{\mathop{\rm Re}\nolimits}

\def\Gr{\mathop{\rm Gr}\nolimits}

\def\Isom{\mathop{\rm Isom}\nolimits}
\def\Diff{\mathop{\rm Diff}\nolimits}
\def\O{\mathop{\rm O}\nolimits}
\def\SO{\mathop{\rm SO}\nolimits}
\def\U{\mathop{\rm U}\nolimits}
\def\SU{\mathop{\rm SU}\nolimits}
\def\PSO{\mathop{\rm PSO}\nolimits}
\def\PSU{\mathop{\rm PSU}\nolimits}
\def\Sp{\mathop{\rm Sp}\nolimits}
\def\Bihol{\mathop{\rm Bihol}\nolimits}

\def\Re{\mathop{\rm Re}\nolimits}
\def\Im{\mathop{\rm Im}\nolimits}

\newcommand\co{\colon} %macro to use in f\co \rightarrow Y
 \def\HollowBoxx #1#2#3{{\dimen0=#1 \advance\dimen0 by -#2
       \dimen1=#1 \advance\dimen1 by #3
        \vrule height 0pt depth #3 width #2
       \hskip -#3
       \vrule height #1 depth #3 width #3}}
 \def\LeftContraction{\mathord{\kern1.45pt \HollowBoxx{6pt}{3.5pt}{.4pt}}\,}

 \def\HollowBox #1#2#3{{\dimen0=#1 \advance\dimen0 by -#3
       \dimen1=#1 \advance\dimen1 by #3
        \vrule height #1 depth #3 width #3
        \vrule height 0pt depth #3 width #2
        \hskip -#3}}
 \def\RightContraction{\mathord{\, \HollowBox{6pt}{3.1pt}{.4pt}} \kern1.6pt}

\newtheorem{theorem}{THEOREM}[section]
\newtheorem{corollary}[theorem]{Corollary}

\newtheorem{proposition}[theorem]{Proposition}

\theoremstyle{definition}

\newtheorem{lemma}[theorem]{Lemma}

\theoremstyle{remark}
\newtheorem{remark}[theorem]{Remark}

\makeatletter
\def\blfootnote{\xdef\@thefnmark{}\@footnotetext}
\makeatother

\begin{document}

\title[Proper group actions in complex geometry]{Proper group actions in
\vspace{0.3cm}\\
complex geometry}
\author[Isaev]{A. V. Isaev}

\address{Mathematical Sciences Institute\\
Australian National University\\
Acton, ACT 2601, Australia}
\email{alexander.isaev@anu.edu.au}

\maketitle

\thispagestyle{empty}

\pagestyle{myheadings}

\begin{abstract} Proper group actions are ubiquitous in mathematics and have many of the attractive features of actions of compact groups. In this survey, we discuss proper actions of Lie groups on smooth manifolds. If the group dimension is sufficiently high, all proper effective actions can be explicitly determined, and our principal goal is to provide a comprehensive exposition of known classification results in the complex setting. They include a complete description of Kobayashi-hyperbolic manifolds with high-dimensional automorphism group, which is a case of special interest. 
\end{abstract}

\tableofcontents

\section{Introduction}\label{intro}
\setcounter{equation}{0}

The purpose of this survey is to give a comprehensive overview of classification results for proper Lie group actions on smooth manifolds, with emphasis on the complex-geometric setting. Given a Lie group $G$ acting on a smooth manifold $X$ by diffeomorphisms, the action is called proper if the map
$$
G\times X\to X\times X,\quad (g,x)\mapsto (gx,x)
$$
is proper. This requirement is easily seen to be satisfied if $G$ is compact, and it is the compact groups whose actions had been traditionally best-understood. However, the result of Palais on the existence of slices (see \cite{Pal}) drew significant attention to the much larger class of proper actions, which turned out to possess many of the good properties of actions of compact groups.

Together with classical results of Myers and Steenrod on groups of Riemannian isometries (see \cite{MySt}), the theorem of Palais implies that the groups acting properly and effectively on a given manifold $X$ can be characterized exactly as the closed subgroups of the groups of isometries with respect to all possible Riemannian metrics on $X$. This fact is of fundamental importance as it relates the theory of proper actions to the classical study of groups of Riemannian motions, which goes back to the foundational work by Ricci, Bianchi, Fubini, \'E.~Cartan and many others. At the same time, one should note that the investigation of groups of motions was primarily conducted in the local setting, whereas proper group actions concern Riemannian isometry groups considered globally.

Analogously to the local setting, the dimension of a Lie group $G$ acting properly and effectively on a smooth $n$-dimensional manifold $X$ does not exceed $n(n+1)/2$. Furthermore, assuming that both $X$ and $G$ are connected, one can show that  for $n\ge 2$ this upper bound is attained only for $X$ isometric to one of the following spaces of constant sectional curvature: the Euclidean space $\RR^n$, the sphere $S^n$, the hyperbolic space $\RR\HH^n$ and the projective space $\RR\PP^n$, with the group $G$ being isomorphic to the connected component of the isometry group of the corresponding model. As a natural extension of this result, one would like to attempt to classify all pairs $(X,G)$ for which $\dim G$ is lower than $n(n+1)/2$ but is still \lq\lq sufficiently high\rq\rq. 

This problem was addressed in a large number of papers spanning a long period from at least as early as the 1920s to the 1970s. Although the term \lq\lq proper actions\rq\rq\,\hspace{-0.1cm} was almost never explicitly mentioned, this is what the authors effectively studied as many results were applicable not only to the local setting, but also to closed subgroups of Riemannian isometry groups treated globally. The most impressive contribution to the area was made by Kobayashi and Nagano who produced, among other results, a classification of pairs $(X,G)$ for group dimension\linebreak $(n-1)(n-2)/2+2$, with $n\ge 6$ (see \cite{KNa}). To this day, it remains to be the lowest value of $\dim G$ for which a complete explicit classification has been found. Beyond this value, the problem of obtaining an effective description of pairs $(X,G)$ is extremely complicated and a solution appears to be out of reach.

In this survey, we focus on classifying pairs $(X,G)$ in the complex setting, i.e., when $X$ is a complex manifold and $G$ lies in the group $\Bihol(X)$ of biholomorphic automorphisms of $X$. Letting $N$ be the complex dimension of $X$ (so that $n=2N$), it is not hard to show that $\dim G\le N^2+2N$. Furthermore, analogously to the smooth case, this upper bound is attained only for $X$ biholomorphically isometric to one of the complete simply-connected K\"ahler spaces of constant holomorphic sectional curvature: the Euclidean space $\CC^N$, the projective space $\CC\PP^N$ and the hyperbolic space $\CC\HH^N$, with the group $G$ being isomorphic to the group of biholomorphic isometries of the corresponding model. As before, one wishes to extend the above result by classifying pairs $(X,G)$ with $\dim G<N^2+2N$ whenever possible. 

Regarding this complex variant of the classification problem, we first of all note that 
$$
N^2+2N<\frac{(n-1)(n-2)}{2}+2\quad \hbox{for $N\ge 5$},
$$
i.e., for $N\ge 5$ the group dimension range relevant to the complex case lies strictly below the dimension range investigated in the smooth situation. Therefore, one cannot rely on the classifications of pairs $(X,G)$ obtained in the smooth case to derive classifications in the complex setting simply by identifying $G$-invariant complex structures on $X$. To deal with the lower values of $\dim G$, one needs to develop new ideas based on the fact that $G$ operates on $X$ by biholomorphisms.                                                                                               

There has been substantial progress in this area, with all pairs $(X,G)$ having been explicitly determined for $\dim G\ge N^2+1$ (see, e.g., \cite{IKru3}). The classification is rather involved despite the fact that in this case the manifold $X$ is automatically $G$-homogeneous. For $\dim G<N^2+1$ the $G$-action need not be transitive, which makes the classification problem even harder. Nevertheless, there are explicit descriptions of manifolds $X$ with $\dim\Bihol(X)=N^2$ and $\dim\Bihol(X)=N^2-1$ under the assumption that $X$ is Kobayashi-hyperbolic, which is a condition guaranteeing that $\Bihol(X)$ is a Lie group acting on $X$ properly. The class of Kobayashi-hyperbolic manifolds is rather important in complex geometry, and finding all manifolds with automorphism group of prescribed non-maximal dimension in this class was in fact our original motivation for studying general proper actions. The hardest part of our classification in this case is describing 2-dimensional manifolds with 3-dimensional automorphism group (see \cite{Isa7}).

The paper is organized as follows. In Sect.~\ref{setup} we review fundamental facts on proper actions starting with the topological setup and gradually progressing to actions of Lie groups on smooth manifolds by diffeomorphisms. The latter are discussed in Sect.~\ref{riemanniancase}, where a number of classification results for pairs $(X,G)$ are presented beginning with the maximal group dimension $n(n+1)/2$. Next, in Sect.~\ref{complexcase} we switch to the complex setting and explain what is known in this case culminating in a detailed classification for $\dim G=N^2+1$. Finally, in 
Sect.~\ref{kobayashihyperboliccase} we specialize to Kobayashi-hyperbolic manifolds and give classifications for $\dim\Bihol(X)=N^2$ and $\dim\Bihol(X)=N^2-1$, with the main result being a complete list of 2-dimensional manifolds whose automorphism groups have dimension 3.

{\bf Acknowledgements.}  This work is supported by the Australian Research Council.

\section{Fundamentals of proper actions}\label{setup}
\setcounter{equation}{0}

Let $X$ be a topological space and ${\mathcal H}(X)$ the group of homeomorphisms of $X$. We endow ${\mathcal H}(X)$ with the {\it compact-open topology}, which is the topology generated by all sets of the form
$$
\left\{f\in{\mathcal H}(X): f(C)\subset U\right\},
$$
where $C\subset X$ is compact and $U\subset X$ is open. For a subgroup\linebreak $H\subset{\mathcal H}(X)$ the subspace topology on $H$ is called the {\it compact-open topology on $H$}. By choosing in the above definition the compact subsets $C$ to be single points, one obtains the {\it topology of pointwise convergence}. Everywhere below, unless stated otherwise, subgroups of ${\mathcal H}(X)$ are considered with the compact-open topology.

Further, a group $G$ equipped with a topology is called a {\it topological group}\, if the maps
\begin{equation}
\begin{array}{ll}
G\times G\to G, & (f,g)\mapsto f\cdot g,\\
\vspace{-0.1cm}\\
G\to G,& f\mapsto f^{-1}
\end{array}\label{topgroup}
\end{equation}
are continuous. For example, if $X$ is locally connected (i.e., has a basis of connected neighborhoods at every point), locally compact and Hausdorff, then ${\mathcal H}(X)$ is a topological group (see \cite{Ar}).

A topological group $G$ is said to {\it act}\, on $X$ if one specifies a continuous map, called an {\it action map}, or simply an {\it action}, 
$$
\Phi\co G\times X\to X
$$
such that for every $g\in G$ one has $\Phi(g,\cdot)\in {\mathcal H}(X)$ and
$$
\Phi(g_1g_2,x)=\Phi(g_1,\Phi(g_2,x))\,\,\hbox{for all $g_1,g_2\in G$, $x\in X$.}
$$
It is customary to abbreviate $\Phi(g,x)$ as $gx$, and using this notation, for $x\in X$ we define the {\it stabilizer}, or {\it isotropy subgroup}, of $x$ as
$$
G_x:=\{g\in G: gx=x\},
$$
and the {\it orbit}\, of $x$ as
$$
Gx:=\{gx\in X:g\in G\}.
$$
An action\ with a single orbit is called {\it transitive}.

Clearly, the action map induces a homomorphism
$$
\hat\Phi:G\to{\mathcal H}(X),\quad g\mapsto \Phi(g,\cdot),
$$
and the action is called {\it effective}\, if $\hat\Phi$ is injective. The continuity of $\Phi$ implies that $\hat\Phi$ is continuous; furthermore, if $X$ is locally compact and Hausdorff, one can equivalently define an action of $G$ on $X$ by specifying a continuous homomorphism $G\to{\mathcal H}(X)$ (see, e.g., \cite{Du}, p.~261 and \cite{Si}, Theorems 11.2.10, 11.2.11). 

In what follows we will be mostly concerned with the situation where $X$ is a smooth manifold (here and below all smooth objects are $C^{\infty}$-smooth). In this case, $G$ is said to {\it act on $X$ by diffeomorphisms}\, if $\hat\Phi(G)$ lies in $\Diff(X)$, the subgroup of ${\mathcal H}(X)$ that consists of all diffeomorphisms of $X$. By the result of \cite{Ar} quoted above, every subgroup of $\Diff(X)$ is a topological group acting on $X$ effectively by diffeomorphisms. Assuming, furthermore, that $X$ is a complex manifold, we say that $G$ {\it acts on $X$ by biholomorphisms}\, if $\hat\Phi(G)$ lies in the group $\Bihol(X)$ of all biholomorphic automorphisms of $X$ (sometimes we refer to this group simply as the {\it automorphism group}). Analogously, every subgroup of $\Bihol(X)$ is a topological group that acts on $X$ effectively by biholomorphisms.

An action of $G$ on $X$ is called {\it proper}\, if the mapping
$$
\tilde\Phi\co G\times X\to X\times X,\quad (g,x)\mapsto (gx,x)
$$
is proper. Recall that a continuous map $f:Y\to Z$ of topological spaces is proper if $f$ is closed and $f^{-1}(z)$ is compact in $Y$ for every\linebreak $z\in Z$. If $Z$ is locally compact and Hausdorff, this condition is equivalent to the requirement that the preimage of every compact subset of $Z$ under $f$ be compact in $Y$ (see, e.g., \cite{Bou}, Chap.~I, pp.~72, 75, 77). 

We will now collect several basic facts on proper actions. Detailed proofs can be found in \cite{Bou}, Chap.~I, \S 10 and Chap.~III, \S 4; \cite{Bil}; \cite{tDie}, pp.~27--29; \cite{L}, pp.~118--121, 318--322.
\vspace{-0.1cm}\\ 

\noindent (I) If $X$ is Hausdorff, then any action of a compact group on $X$ is proper.
\vspace{-0.1cm}\\

\noindent (II) If $G$ acts properly on $X$, then so does any closed subgroup of $G$.
\vspace{-0.1cm}\\  

\noindent (III) Considering the preimage under $\tilde\Phi$ of a point $\{(x,x)\}\subset X\times X$, we see that the properness of an action of $G$ on $X$ implies that the stabilizer $G_x$ of $x$ is compact in $G$. Notice that the compactness of all the stabilizers is not in general equivalent to properness. In fact, the following holds:

\begin{proposition}\label{propgentopspace} Let a topological group $G$ act on a Hausdorff topological space $X$. Then the action is proper if and only if for every compact subset $C\subset X$ the set $\{g\in G: gC\cap C \ne \emptyset\}$ is compact in $G$. Nevertheless, for any compact subgroup $K\subset G$, the natural action of $G$ on the quotient $G/K$ is proper.
\end{proposition}

\noindent (IV) If $G$ acts properly on $X$, then for every $x\in X$ the map
$$
\Phi(\cdot,x)\co G\to X, \quad g\mapsto gx
$$
is proper. This yields: 
\vspace{-0.3cm}\\

\hspace{0.5cm} (IVa) the orbit $Gx$ of $x$ is closed in $X$, and the natural bijection $G/G_x\to Gx$ is a homeomorphism;
\vspace{-0.3cm}\\

\hspace{0.5cm} (IVb) if $X$ is locally compact, then so is $G$, and $\hat\Phi(G)$ is closed\linebreak in ${\mathcal H}(X)$.
\vspace{-0.1cm}\\

\noindent (V) Let $G$ be Hausdorff and
$$
H:=\bigcap_{x\in G}G_x
$$
be the {\it kernel}\, of a proper action of $G$ on $X$. Clearly, $H$ is a normal closed subgroup of $G$, and the Hausdorff topological group $G/H$ acts on $X$ properly and effectively. Then the continuous homomorphism 
$$
\hat\Phi_H:G/H\to{\mathcal H}(X)
$$
induced by this action is a homeomorphism from $G/H$ onto its image\linebreak $\hat\Phi_H(G/H)=\hat\Phi(G)$. In addition, the compact-open topology and the topology of pointwise convergence on $\hat\Phi(G)$ coincide.
\vspace{-0.1cm}\\

Perhaps the first example of a large class of not necessarily compact groups acting properly was given by the following result (see \cite{vDvdW} and \cite{KNo}, Chap.~I, Theorem 4.7):

\begin{theorem}\label{vdvdw} If $(X,{\mathtt d})$ is a connected locally compact metric space, then the group $\Isom(X,{\mathtt d})\subset{\mathcal H}(X)$ of isometries of $X$ with respect to the distance ${\mathtt d}$ is a topological group acting properly on $X$.
\end{theorem}

We now specialize to the case when $X$ is a connected smooth manifold and $G$ a topological group acting on $X$ properly and effectively by diffeomorphisms. Note that $G$ is Hausdorff  since ${\mathcal H}(X)$ is. Our broad aim is to understand such actions. We start with the following important fact:

\begin{theorem}\label{Liegroup} Let $H$ be a locally compact group acting effectively by diffeomorphisms on a connected smooth manifold $Y$. Then $H$ has the structure of a Lie group.
\end{theorem}

By a Lie group in the above statement we mean a smooth manifold (not necessarily second countable) such that the operations of multiplication  and taking the inverse are smooth maps (cf.~(\ref{topgroup})). In particular, a Lie group so defined may have uncountably many connected components. Nevertheless, any Lie group is paracompact, i.e., its every connected component is second countable (see, e.g., \cite{H}, p.~88).

Theorem \ref{Liegroup} was first obtained in \cite{BMo2} and \cite{Kur} under a hypothesis more restrictive than the effectiveness of the action. Namely, it was required that if for $h\in H$ its fixed-point set $\{y\in Y:hy=y\}$ has an interior point, then $h=e$, the identity element of $H$. In fact, Theorem 1 in \cite{BMo2} shows that if $H$ acts effectively, then it {\it has no small subgroups}, which means that there exists a neighborhood of $e$ in $H$ containing no non-trivial subgroups (see also \cite{MZ2}, Sect.~5.2). Therefore, Theorem \ref{Liegroup} in full generality is a consequence of Theorem 1 in \cite{BMo2} and later results obtained in articles \cite{G} and \cite{Yam}, where the condition of having no small subgroups was utilized in an essential way (see \cite{T}, Chap.~5 for a review of these ideas). [We note that article \cite{G} was instrumental in solving Hilbert's fifth problem in \cite{MZ1} and that Theorem \ref{Liegroup} confirms the Hilbert-Smith conjecture for actions by diffeomorphisms on smooth manifolds.] Further, by \cite{BMo1} the group $H$ {\it acts on $Y$ smoothly}, i.e., the corresponding action map is smooth (see also \cite{MZ2}, Sect.~5.2). 

Combining properties (IV), (V) with Theorem \ref{Liegroup} and the result of \cite{BMo1}, we can assume that:

\begin{itemize}

\item[(VI)] as an abstract group, $G$ is a subgroup of $\Diff(X)\subset{\mathcal H}(X)$, thus $\hat\Phi$ is simply the inclusion map;

\item[(VII)] the topology of $G$ is the compact-open topology (in particular, $G$ is second countable when $X$ is paracompact), which coincides with the pointwise topology as well as with that induced by the weak Whitney topology of uniform convergence of diffeomorphisms and all their partial derivatives on compact subsets\linebreak of $X$;

\item[(VIII)] $G$ is closed in ${\mathcal H}(X)$;

\item[(IX)] $G$ is a Lie group, which has no more than countably many connected components when $X$ is paracompact;

\item[(X)] the action of $G$ on $X$ is smooth;

\item[(XI)] the $G$-orbits are closed smooth submanifolds of $X$. 

\end{itemize}

From now on we additionally suppose that $X$ is paracompact. It turns out that in this case proper actions can be described in the language of Riemannian geometry. In what follows we denote a smooth manifold $Y$ endowed with a smooth Riemannian metric ${\mathtt f}$ by $[Y,{\mathtt f}]$.  Recall that any Riemannian manifold is automatically paracompact (see, e.g. \cite{H}, p.~52). 

We need the following result (see \cite{Pal}, Theorem 4.3.1; \cite{Al}; \cite{DK}, Proposition 2.5.2):

\begin{theorem}\label{invRiemmetric} If a Lie group $H$ acts smoothly and properly on a smooth paracompact manifold $Y$, then $Y$ admits a smooth $H$-invariant Riemannian metric. Furthermore, if $Y$ is complex and $H$ acts by biholomorphisms, then one can choose the metric to be real-analytic and Hermitian. 
 \end{theorem}
 
\noindent Taking into account the assumptions on the group $G$ summarized above, we see that by Theorem \ref{invRiemmetric} there exists a smooth Riemannian metric ${\mathtt h}$ on $X$ such that $G$ is a closed subgroup of the group $\Isom[X,{\mathtt h}]\subset\Diff(X)$ of {\it Riemannian isometries}\, of $X$, i.e., diffeomorphisms of $X$ preserving ${\mathtt h}$. We will now explain that this is in fact a complete description of all groups acting properly on $X$.

Suppose that $X$ is endowed with a smooth Riemannian metric ${\mathtt g}$. One can then consider the distance associated to ${\mathtt g}$ as follows:
$$
{\mathtt d}_{\mathtt g}(x,y):=\inf_{\gamma}\int_0^1||\gamma'(t)||_{\mathtt g}\,dt,\quad x,y\in X,
$$
where the inf is taken over all piecewise regular curves $\gamma:[0,1]\ra X$, such that $\gamma(0)=x$, $\gamma(1)=y$. It is well-known (see, e.g., \cite{H}, Chap.~I, \S 9) that ${\mathtt d}_{\mathtt g}$ is indeed a distance and induces the topology of $X$.\footnote{Note that in general any continuous inner distance on a locally compact Hausdorff space induces its topology (see \cite{Ba2}).} The distance ${\mathtt d}_{\mathtt g}$ turns $X$ into a connected locally compact metric space, and therefore Theorem \ref{vdvdw} implies that $\Isom(X,{\mathtt d}_{\mathtt g})$ is a topological group acting properly on $X$. 

Further, observe that $\Isom[X,{\mathtt g}]\subset\Isom(X,{\mathtt d}_{\mathtt g})$. Remarkably, the following holds (see \cite{MySt} and \cite{H}, Chap.~I, \S11):

\begin{theorem}\label{myerssteenrod} For a connected Riemannian manifold $[Y,{\mathtt f}]$ one has
$$  
\Isom[Y,{\mathtt f}]=\Isom(Y,{\mathtt d}_{\mathtt f}).
$$
\end{theorem}      

\noindent By Theorem \ref{myerssteenrod}, we have $\Isom[X,{\mathtt g}]=\Isom(X,{\mathtt d}_{\mathtt g})$, hence $\Isom[X,{\mathtt g}]$ is a topological group acting properly on $X$. Therefore, the action of every closed subgroup of $\Isom[X,{\mathtt g}]$ on $X$ is proper as well.

We can now summarize our discussion as follows: describing all proper effective actions of topological groups by diffeomorphisms on a connected smooth paracompact manifold $X$ means describing all closed subgroups of the groups of isometries $\Isom[X,{\mathtt g}]$, where ${\mathtt g}$ runs over all smooth Riemannian metrics on $X$ (cf.~\cite{DK}, p.~106).    

\section{Proper actions of high-dimensional groups}\label{riemanniancase}
\setcounter{equation}{0}

Our broad aim is to classify, whenever possible, pairs $(X,G)$, where $X$ is a connected smooth paracompact manifold and $G\subset\Diff(X)$ a group acting properly, with two such pairs $(X_1,G_1)$ and $(X_2,G_2)$ called {\it equivalent}\, if there exists a diffeomorphism $f:X_1\to X_2$ satisfying $f\circ G_1\circ f^{-1}=G_2$. In what follows we always assume $G$ to be connected. This assumption is convenient for technical reasons and does not significantly impact on the generality of the results. For any Lie group $H$, we denote by $H^{\circ}$ its connected component containing the identity, thus for a $G$-invariant Riemannian metric ${\mathtt g}$ on $X$ one has $G\subset\Isom[X,{\mathtt g}]^{\circ}$.

It is natural to start with \lq\lq large\rq\rq\hspace{-0.01cm} groups $G$ as manifolds having \lq\lq many\rq\rq\hspace{-0.01cm} symmetries are easier to describe. There are various ways in which one can specify what being \lq\lq large\rq\rq\hspace{-0.01cm} exactly means. Our approach is to assess the \lq\lq size\rq\rq\hspace{-0.01cm} of $G$ in terms of its dimension $\dim G$ as a Lie group. Setting $n:=\dim X$, we will now estimate $\dim G$ in terms of $n$.

First, we state the important fact known as Bochner's linearization theorem, which will be a useful tool throughout the entire survey (see \cite{Boc}; \cite{DK}, Theorem 2.2.1; \cite{MZ2}, pp.~206--208). 

\begin{theorem}\label{Bochner'slineartheorem} Let $Y$ be a smooth manifold and $K\subset\Diff(Y)$ a compact subgroup fixing a point $y\in Y$. Then there exists a $K$-invariant neighborhood $U$ of $y$ in $Y$, a neighborhood $V$ of the origin in the tangent space $T_y(Y)$ to $Y$ at $y$ and a diffeomorphism $\varphi: U\to V$ such that for every $h\in K$ one has
$$
d_yh\circ \varphi=\varphi\circ h|_U.
$$
\end{theorem}   

\noindent In other words, Theorem \ref{Bochner'slineartheorem} asserts that near the fixed point $y$ the action of the group $K$ can be identified with the action of the subgroup of $\GL(T_y(Y),\RR)$ that consists of the linear parts of the elements of $K$ at $y$.

\begin{remark}\label{bochnercomplex} If in Theorem \ref{Bochner'slineartheorem} the manifold $Y$ is complex and\linebreak $K\subset\Bihol(Y)$, then $K$ is identified with a subgroup of $\GL(T_y(Y),\CC)$. Furthermore, the proof of the theorem, as presented, for example, in \cite{DK}, yields that in this case one can choose $\varphi: U\to V$ to be a biholomorphism (see \cite{BMa}, p.~55). For bounded domains in $\CC^N$ this complex version of the theorem is related to classical work of H. Cartan (see \cite{CaH}, Art.~13, Chap.~II, IV and Art.~14; \cite{Mar}). Also, as shown in \cite{Ka}, the proof can be generalized to complex spaces (cf.~\cite{Ak},\linebreak pp.~36--37).  
\end{remark} 

Fix $x\in X$ and consider the {\it isotropy representation}\, $\alpha_x$ of the stabilizer $G_x$:
$$
\alpha_x: G_x\to \GL(T_x(X),\RR),\quad g\mapsto d_xg.
$$
Properties (I) and (VII) in Sect.~\ref{setup} specialized to the action of the compact group $G_x$ imply that $\alpha_x$ is continuous. Furthermore, application of Theorem \ref{Bochner'slineartheorem} to $G_x$ yields that $\alpha_x$ is faithful. It then follows that the isotropy representation is a homeomorphism onto its image. In particular, the {\it linear isotropy subgroup}\, $LG_x:=\alpha_x(G_x)$ is isomorphic to $G_x$ as a Lie group. 

Since the subgroup $LG_x\subset \GL(T_x(X),\RR)$ is compact, one can find a scalar product on $T_x(X)$ preserved by $LG_x$ (moreover, by Theorem \ref{invRiemmetric} there exists a $G$-invariant Riemannian metric on $X$). Choosing coordinates in $T_x(X)$ in which the scalar product is given by the identity matrix, we see that $LG_x$ can be realized as a closed subgroup of the orthogonal group $\O_n(\RR)$. This implies
\begin{equation}
\dim G_x=\dim LG_x\le \dim \O_n(\RR)=\frac{n(n-1)}{2}.\label{dimstab}
\end{equation}
At the same time, we have
\begin{equation}
\dim Gx\le \dim X=n.\label{dimorbit}
\end{equation}
From inequalities (\ref{dimstab}) and (\ref{dimorbit}) we derive the following estimate for $\dim G$:
$$
\dim G=\dim Gx+\dim G_x\le \frac{n(n+1)}{2}.
$$

Note that classifying pairs $(X,G)$ for $n=1$ is easy. Indeed, in this case the group $G$ is either one-dimensional or trivial, with $X$ being diffeomorphic to one of $\RR$, $S^1$. If $\dim G=1$, then $(X,G)$ is equivalent either to $(\RR,\RR)$ (where $\RR$ acts on itself by translations) or to $(S^1,\SO_2(\RR))$ (where $\SO_2(\RR)$ acts on $S^1$ by rotations). Thus, from now on we assume that $n\ge 2$. 

For the maximal group dimension $n(n+1)/2$, one has the following classification (see \cite{KNo}, p.~308; \cite{Ko1}, Chap.~II, Sect.~3; \cite{O}):

\begin{theorem}\label{maxgroupdim} If $\dim G=n(n+1)/2$, then the pair $(X,G)$ is equivalent to one of the pairs
$$
\begin{array}{cl}
\hbox{\rm (i)}& (\RR^n,\SO_n(\RR)\ltimes \RR^n),\\
\vspace{-0.3cm}\\
\hbox{\rm (ii)} & (S^n,\SO_{n+1}(\RR)),\\
\vspace{-0.3cm}\\
\hbox{\rm (iii)} & (\RR\HH^n, \SO_{n,1}^{\circ}),\\
\vspace{-0.3cm}\\
\hbox{\rm (iv)} & (\RR\PP^n,\PSO_{n+1}(\RR)),
\end{array}
$$
where $\SO_n(\RR)\ltimes \RR^n$ is the group of orientation-preserving motions of the Euclidean space $\RR^n$, the groups $\SO_{n+1}(\RR)$ and $\PSO_{n+1}(\RR)$ are regarded as subgroups of the groups of isometries of the sphere $S^n$ and projective space $\RR\PP^n$, respectively, and $\SO_{n,1}^{\circ}$ is regarded as a subgroup of the group of isometries of the hyperbolic space $\RR\HH^n$. 
\end{theorem}

We will now briefly outline the proof of Theorem \ref{maxgroupdim}. Since the dimension of $G$ is maximal possible, its action on $X$ is transitive and for every $x\in X$ the group $LG_x$ contains an open subgroup that in some coordinates in $T_x(X)$ coincides with $\SO_n(\RR)$.  Notice that $\SO_n(\RR)$ acts transitively on the Grassmannian $\Gr(2,n)$ of 2-planes in $\RR^n$. Fixing a $G$-invariant Riemannian metric ${\mathtt g}$ on $X$, we then deduce that the sectional curvature of ${\mathtt g}$ is constant. Since $[X,{\mathtt g}]$ is {\it Riemannian homogeneous}\, (i.e., $\Isom[X,{\mathtt g}]$ acts on $X$ transitively), it follows that $[X,{\mathtt g}]$ is {\it complete}, i.e., the metric space $(X,{\mathtt d}_{{\mathtt g}})$ is complete (see \cite{KNo}, Chap.~IV, Theorem 4.5). Therefore, if $X$ is simply-connected, $[X,{\mathtt g}]$ is isometric to one of $\RR^n$, $S^n$, $\RR\HH^n$ endowed with a metric proportional to the Euclidean, spherical, or hyperbolic one, respectively (see \cite{KNo}, Chap.~VI, Theorem 7.10). If $X$ is not simply-connected, one considers its universal cover $\tilde X$, and it is not hard to observe that $[\tilde X,\tilde{\mathtt g}]$ (with $\tilde{\mathtt g}$ being the pull-back of ${\mathtt g}$) must be isometric to $S^n$, which leads to $[X,{\mathtt g}]$ being isometric to $\RR\PP^n$. Since $G=\Isom[X,{\mathtt g}]^{\circ}$, the isometry so constructed transforms the group $G$ into the identity component of the isometry group of one of $\RR^n$, $S^n$, $\RR\HH^n$, $\RR\PP^n$, and we obtain the classification stated in Theorem \ref{maxgroupdim}. Note that the constancy of the sectional curvature of a Riemannian manifold whose isometry group has maximal possible dimension has been known for a long time (see, e.g., \cite{CaE}, p.~269; \cite{Ei}, p.~216; \cite{Yan2}, Theorem 6.2).

\begin{remark}\label{anotherproof} In the above proof, since $LG_x$ acts transitively on the unit sphere in $T_x(X)$ for all $x\in X$, the Riemannian manifold $[X,{\mathtt g}]$ is {\it isotropic} (i.e., the group of the linear parts of the isometries fixing $x$ acts transitively on the unit sphere in $T_x(X)$). Therefore, Theorem \ref{maxgroupdim} can be also derived from the classification of isotropic Riemannian manifolds (see \cite{Wo}, Theorems 8.12.2 and 8.12.8).
\end{remark}

An important observation in the proof of Theorem \ref{maxgroupdim} is the fact that for all $x\in X$ the linear isotropy subgroup $LG_x$ contains a copy of $\SO_n(\RR)$. In general, the \lq\lq size\rq\rq\hspace{-0.01cm} and explicit form of the linear isotropy subgroups are rather essential for obtaining classifications of pairs\linebreak $(X,G)$ and will be a recurrent theme throughout the survey. 

Clearly, for any $x\in X$ we have
\begin{equation}
\dim LG_x=\dim G-\dim Gx\ge \dim G-n,\label{stabestim}
\end{equation}
and therefore, by imposing a lower bound on $\dim G$, one can control $\dim LG_x$ from below. For example, assuming 
\begin{equation}
\dim G\ge\frac{n(n-1)}{2}+2,\quad n\ge 3,\label{estimG}
\end{equation}
from (\ref{stabestim}) we obtain
\begin{equation}
\dim LG_x\ge\frac{(n-1)(n-2)}{2}+1\label{estimLG}
\end{equation}
for any $x\in X$. Interestingly, one has (see \cite{MoSa}):

\begin{lemma}\label{monsamellemma} \it For $n\ge 2$, $n\ne 4$ the only proper closed subgroup of $\O_n(\RR)$ of dimension at least $(n-1)(n-2)/2+1$ is $\SO_n(\RR)$. 
\end{lemma}

Using Lemma \ref{monsamellemma}, one can easily prove that if $G$ satisfies (\ref{estimG}), it acts transitively on $X$ (see \cite{Ko1}, p.~48), and therefore we have\linebreak $\dim G=n(n+1)/2$, which is the case already considered in Theorem \ref{maxgroupdim}. In particular, the following holds (see \cite{Wan}, \cite{Yan1}, \cite{Eg1}):

\begin{theorem}\label{firstlacuna} For $n\ge 3$, $n\ne 4$ a group $G$ with
$$
\frac{n(n-1)}{2}+2\le\dim G<\frac{n(n+1)}{2}
$$
cannot act properly on $X$.
\end{theorem}

Theorem \ref{firstlacuna} shows that for $n\ge  3$, $n\ne 4$ the set of the dimensions of the groups that can act properly on $X$ has a lacuna of size linear in $n$ located immediately below the maximal possible dimension. Such lacunas are common in classification problems of this kind (see, e.g., \cite{Eg3}; \cite{Eg4}, p.~ 219; \cite{KNa}; \cite{Man}; \cite{Wak}), and we will encounter a similar gap phenomenon in the complex setting considered in the next section (see Corollary \ref{firstlacunacomplexcor}). We also note that Theorem \ref{firstlacuna} extends -- in the global setting -- a result of \cite{F1} (see p.~54 therein) on the non-existence of actions of groups of dimension $n(n+1)/2-1$ for $n\ge 3$.

For $n=4$ inequalities (\ref{estimG}) and (\ref{estimLG}) become $\dim G\ge 8$ and $\dim LG_x\ge 4$, respectively. Although the group $\SO_4(\RR)$ is excluded from Lemma \ref{monsamellemma}, there is a complete description of all connected Lie subgroups of $\SO_4(\RR)$ (see \cite{Ish}). It shows, in particular, that $\SO_4(\RR)$ has no proper subgroups of dimension at least 5 and all subgroups of dimension 4 are conjugate in $\O_4(\RR)$ to the unitary group $\U_2\subset\SO_4(\RR)$. This rules out the case $\dim G=9$ and leads to a complete classification for $\dim G=8$ as follows (see \cite{Ish}, \cite{Eg2} and Theorem \ref {maxdimcompl} in the next section): 

\begin{theorem}\label{classforn=4} If $n=4$, then $\dim G\ne 9$, and if $\dim G=8$, then $X$ has a $G$-invariant complex structure. In the latter case, the pair $(X,G)$ is equivalent {\rm (}by means of a biholomorphism{\rm )} to one of the pairs
$$
\begin{array}{cl}
\hbox{\rm (i)}& (\CC^2,\U_2\ltimes \CC^2),\\
\vspace{-0.3cm}\\
\hbox{\rm (ii)} & (\CC\PP^2,\PSU_3),\\
\vspace{-0.3cm}\\
\hbox{\rm (iii)} & (\BB^2, \PSU_{2,1}),
\end{array}
$$
where $\U_2\ltimes \CC^2$ is a subgroup of the group $\SO_4(\RR)\ltimes\RR^4$ of orientation-preserving motions of the Euclidean space $\RR^4=\CC^2$, the group $\PSU_3$ is regarded as a subgroup of the group of isometries of $\CC\PP^2$ with respect to the Fubini-Study metric, and $\BB^2$ is the unit ball in $\CC^2$ for which the group $\Bihol(\BB^2)$ is identified with $\PSU_{2,1}$.
\end{theorem}

\noindent In relation to the classification for $\dim G=8$ in Theorem \ref{classforn=4}, we note that some of the results of \cite{F2} concerning the groups of motions of 4-dimensional manifolds are not entirely correct (see \cite{Eg2} for details).

\begin{remark}\label{complex=2} Recall (see, e.g., \cite{Ru}, Theorem 2.2.5) that for the unit ball $\BB^N$ in $\CC^N$ the group $\Bihol(\BB^N)$ consists of all transformations of the form 
$$
z\mapsto\displaystyle\frac{Az+b}{cz+d},\quad \hbox{with $z:=\left(\hspace{-0.1cm}\begin{array}{c}z_1\\\vdots\\z_N\end{array}\hspace{-0.1cm}\right)$,}\\
$$
where
$$
\left(\hspace{-0.1cm}\begin{array}{cc}
A& b\\
c& d
\end{array}
\hspace{-0.1cm}\right)
\in \SU_{N,1},
$$
which shows that $\Bihol(\BB^N)$ and $\PSU_{N,1}$ are isomorphic for any $N$. Notice that $\BB^N$ endowed with the Bergman metric is a model of the complex $N$-dimensional hyperbolic space $\CC\HH^N$ and that the Bergman metric is $\PSU_{N,1}$-invariant. Therefore, Theorem \ref{classforn=4} is the exact complex analogue of Theorem \ref{maxgroupdim} for $n=2$. A similar statement in an arbitrary complex dimension appears in Theorem \ref{maxdimcompl} in the next section.
\end{remark}

In view of Theorems \ref{firstlacuna}, \ref{classforn=4}, the next group dimension to be studied for $n\ge 3$ is $n(n-1)/2+1$. In this case, for $n\ge 5$ there also exists a complete classification of pairs $(X,G)$ as stated, for example, in \cite{Ko1}, p.~54 (see \cite{O}, \cite{Kui}, \cite{Yan1}). According to the classification, $X$ is isometric, with respect to a $G$-invariant Riemannian metric, either to $\RR\HH^n$, or to the product of one of $\RR$, $S^1$ and one of $\RR^{n-1}$, $S^{n-1}$, $\RR\HH^{n-1}$, $\RR\PP^{n-1}$, and under the isometry the group $G$ is transformed into one of several explicitly given groups. Since the descriptions of some of the groups are somewhat lengthy, we do not provide a complete statement here and refer the reader to \cite{Ko1}, Chap.~II, Sect.~3 for full details.

Notice that in the discussion following Theorem \ref{maxgroupdim} the case $n=2$ was ignored, and in the classification for $\dim G=n(n-1)/2+1$ outlined above the cases $n=3,4$ were omitted. We will now briefly mention known results in these situations. First, for $n=2$, with the exclusion of the maximal dimension, $\dim G$ can only take the values 0, 1, 2. If $\dim G=2$, it is easy to show that the action of $G$ is transitive (see \cite{Ko1}, p.~49), thus the sectional curvature of $X$ is constant. Also, 2-dimensional Riemannian manifolds with positive-dimensional group of isometries were  to some extent described in \cite{My} and \cite{Eb}, Theorem 5.1. Next, if $n=3$ and $\dim G=4$, then again arguing as in \cite{Ko1}, we see that the action of $G$ is transitive, and for a simply-connected $X$ a classification of pairs $(X,G)$ is then given by a theorem stated in \cite{CaE}, p.~303 (see also \cite{Pat}). Finally, the case when $n=4$ and $X$ is simply-connected was considered in \cite{Ish}, where classifications for\linebreak $\dim G=6$ and $\dim G=7$ were obtained. Here for $\dim G=7$ the action is automatically transitive, whereas for $\dim G=6$ transitivity is an assumption. We also mention that low-dimensional cases in the local setting were extensively studied in many classical papers (see, e.g., \cite{Bia}, \cite{F1}, \cite{F2} and the references provided in \cite{CaE}, p.~296).

Next, there exists an explicit classification of pairs $(X,G)$ with\linebreak $(n-1)(n-2)/2+2\le\dim G\le n(n-1)/2$ for $n\ge 6$ (see \cite{KNa} and \cite{Wak}), which is far too lengthy to be reproduced here. There are also many other remarkable results, especially for actions of compact groups (see, e.g., \cite{Man} and \cite{J}, Chap.~IV), but, to the best of our knowledge, no complete explicit classifications beyond group dimension $(n-1)(n-2)/2+2$, with $n\ge 6$, have been found. In this survey, we are primarily interested in actions by biholomorphisms on complex manifolds, and this is the topic we will focus on in the forthcoming sections. For more detail on proper actions in the smooth case we refer the reader to monographs \cite{Ko1} and \cite{Yan2}.

 \section{The complex case}\label{complexcase}
\setcounter{equation}{0}

From now on we concentrate on studying proper group actions in the complex setting. Everywhere below $X$ denotes a connected complex paracompact manifold of complex dimension $N$ (hence $n=2N$) and $G$ a connected closed subgroup of $\Bihol(X)$ acting on $X$ properly. We will be classifying pairs $(X,G)$ up to the following equivalence: two pairs $(X_1,G_1)$ and $(X_2,G_2)$ are said to be {\it biholomorphically equivalent}\,  if there exists a biholomorphism $f:X_1\to X_2$ with $f\circ G_1\circ f^{-1}=G_2$. All Riemann surfaces admitting actions of positive-dimensional groups are well-known (see, e.g., \cite{FK}, p.~294), thus we will be interested in the case $N\ge 2$, although some of the facts stated below apply in the one-dimensional situation as well.

Proper actions by biholomorphisms are found in abundance, as shown by the following result, which holds even for complex spaces (see \cite{Ka}, Satz 2.5):

\begin{theorem}\label{complexmetrics} Let $Y$ be a connected complex manifold and $H$ a closed subgroup of $\Bihol(Y)$. Assume that $H$ preserves a continuous distance on $Y$. Then the action of $H$ on $Y$ is proper.
\end{theorem}

\noindent Together with Theorem \ref{invRiemmetric}, this fact yields that if $Y$ is paracompact, the groups acting properly and effectively on $Y$ by biholomorphisms are exactly the closed subgroups of $\Bihol(Y)$ that preserve continuous distances on $Y$. In the next section, we will consider the important class of Kobayashi-hyperbolic complex manifolds, which includes, for instance, all bounded domains in $\CC^N$. Every such manifold carries a continuous biholomorphically invariant distance, and therefore if $Y$ is Kobayashi-hyperbolic, then $\Bihol(Y)$ acts on $Y$ properly. 

Since $G$ lies in $\Bihol(X)$, for every $x\in X$ the group $LG_x$ is a compact subgroup of $\GL(T_x(X),\CC)$ (and the action of this subgroup near the origin is identified with that of $G_x$ near $x$ by means of a biholomorphism -- see Remark \ref{bochnercomplex}). Hence, $LG_x$ preserves a Hermitian scalar product on $T_x(X)$ (moreover, by Theorem \ref{invRiemmetric} there exists a $G$-invariant Hermitian metric on $X$). Choosing complex coordinates in $T_x(X)$ in which the scalar product is given by the identity matrix, we see that $LG_x$ is isomorphic to a closed subgroup of the unitary group $\U_N$. We then have
\begin{equation}
\dim G_x=\dim LG_x\le \dim \U_N=N^2.\label{dimstabcomplex}
\end{equation}
Together with (\ref{dimorbit}), inequality (\ref{dimstabcomplex}) yields
\begin{equation}
\dim G=\dim Gx+\dim G_x\le N^2+2N.\label{dimgroupcomplex}
\end{equation}   

It is instructive to compare the upper bound on $\dim G$ in (\ref{dimgroupcomplex}) with the lowest group dimension $(n-1)(n-2)/2+2$ (where $n\ge 6$), for which a classification of proper actions by diffeomorphisms on smooth manifolds is known (see Sect.~\ref{riemanniancase}). We have 
\begin{equation}
N^2+2N<\frac{(n-1)(n-2)}{2}+2\quad \hbox{for $N\ge 5$.}\label{inequalityoverlap}
\end{equation}
Thus, the group dimension range that arises in the complex case, for $N\ge 5$ lies strictly below the dimension range investigated in the classical smooth case and therefore is not in any way covered by the results discussed in Sect.~\ref{riemanniancase}. Furthermore, overlaps with these results for $N=2,3,4$ occur in relatively easy situations and do not lead to any significant simplifications in the complex setting (see Remark \ref{overlapsrealcase}). 

This section contains a complete biholomorphic classification of pairs $(X,G)$ for $N^2+1\le\dim G\le N^2+2N$ obtained independently of the facts presented in Sect.~\ref{riemanniancase}. Notice that by Satz 1.2 in \cite{Ka}, in this situation the action of $G$ on $X$ is transitive. All complex manifolds of dimensions 2 and 3 admitting transitive Lie group actions by biholomorphisms were found in \cite{HL}, \cite{OR}, \cite{Wi}, thus we are mainly concerned with the case $N\ge 4$ (although this constraint is not imposed in the theorems that appear below).

First of all, we state a result for the maximal possible dimension\linebreak $N^2+2N$, which also holds for complex spaces (see \cite{Ka} and \cite{Ak}, pp.~49--50). Taking into account Remark \ref{complex=2}, we see that it is the exact complex analogue of Theorem \ref{maxgroupdim}.

\begin{theorem}\label{maxdimcompl} If $\dim G=N^2+2N$, then the pair $(X,G)$ is biholomorphically equivalent to one of the pairs
$$
\begin{array}{cl}
\hbox{\rm (i)}& (\CC^N,\U_N\ltimes \CC^N),\\
\vspace{-0.3cm}\\
\hbox{\rm (ii)} & (\CC\PP^N,\PSU_{N+1}),\\
\vspace{-0.3cm}\\
\hbox{\rm (iii)} & (\BB^N, \PSU_{N,1}),
\end{array}
$$
where $\U_N\ltimes \CC^N$ is a subgroup of the group $\SO_{2N}(\RR)\ltimes\RR^{2N}$ of orientation-preserving motions of the Euclidean space $\RR^{2N}=\CC^N$, the group $\PSU_{N+1}$ is regarded as a subgroup of the group of isometries of $\CC\PP^N$ with respect to the Fubini-Study metric, and $\BB^N$ is the unit ball in $\CC^N$ for which the group $\Bihol(\BB^N)$ is identified with $\PSU_{N,1}$.
\end{theorem}

We will now sketch the proof of Theorem \ref{maxdimcompl}. Fix a $G$-invariant Hermitian metric ${\mathtt g}$ on $X$ and consider the group of {\it biholomorphic isometries}\, of $[X,{\mathtt g}]$ defined as
$$
\Isom_{\Bihol}[X,{\mathtt g}]:=\Isom[X,{\mathtt g}]\cap\Bihol(X).
$$
Clearly, $\Isom_{\Bihol}[X,{\mathtt g}]$ is a closed subgroup of $\Isom[X,{\mathtt g}]$, hence it acts on $X$ properly, which implies $G=\Isom_{\Bihol}[X,{\mathtt g}]^{\circ}$. Further, since the dimension of $G$ is maximal possible, the group $LG_x$ in some complex coordinates in $T_x(X)$ coincides with $\U_N$ for every $x\in X$. Hence, $LG_x$ acts transitively on the unit sphere in $T_x(X)$ for all $x\in X$, and therefore the Riemannian manifold $[X,{\mathtt g}]$ is isotropic. The classification of isotropic Riemannian manifolds (see \cite{Wo}, Theorems 8.12.2 and 8.12.8) now implies that $[X,{\mathtt g}]$ is biholomorphically isometric to one of $\CC^N$, $\CC\PP^N$, $\BB^N$ endowed with a metric proportional to the Euclidean, Fubini-Study, or Bergman metric, respectively. Since $G=\Isom_{\Bihol}[X,{\mathtt g}]^{\circ}$, the biholomorphic isometry so constructed transforms the group $G$ into the (connected) group of biholomorphic isometries of one of $\CC^N$, $\CC\PP^N$, $\BB^N$, and we obtain the classification stated in Theorem \ref{maxdimcompl}.

\begin{remark}\label{otherproofs}
In the above proof we utilized the property that $LG_x$ acts transitively on the unit sphere in $T_x(X)$ for all $x\in X$. In fact, as shown in \cite{GK} and \cite{BDK}, for a connected complex paracompact manifold $Y$ the existence of just one point $y_0\in Y$ fixed by the action of a compact group $K\subset\Bihol(Y)$, with $LK_{y_{{}_0}}$ acting transitively on complex directions in $T_{y_{{}_0}}(Y)$, is sufficient to conclude that $Y$ is biholomorphic to one of $\CC^N$, $\CC\PP^N$, $\BB^N$.
\end{remark}

Further, in order to investigate lower group dimensions, we need the following result, which can be regarded as a unitary analogue of Lemma \ref{monsamellemma}, as well as parts (I), (II) of the lemma in \cite{Ish}, p.~347 (see \cite{IKra2} and \cite{Isa6}, Lemma 1.4):

\begin{lemma}\label{lmlacuna}\it Let $H\subset \U_N$ be a connected closed subgroup, with $N\ge 2$. Suppose that $\dim H\ge N^2-2N+2$. Then either $H=\U_N$, or $H=\SU_N$, or $H$ is conjugate in $\U_N$ to $\U_{N-1}\times \U_1$ {\rm (}realized by matrices in block-diagonal form{\rm )}, or $N=4$ and $H$ is conjugate in $\U_4$ to either $\Sp_2$ or $e^{\hbox{\tiny $i\RR$}}\Sp_2$, where $\Sp_2$ is the compact real form of $\Sp_4(\CC)$.
\end{lemma}

Lemma \ref{lmlacuna} yields a description of the identity components of the linear isotropy subgroups for $\dim G\ge N^2+2$ and, with some additional work, leads to the following classification (see \cite{Isa5}, \cite{Isa8}):

\begin{theorem}\label{firstlacunacomplex} If $N^2+2\le\dim G<N^2+2N$, with $N\ge 2$, then the pair $(X,G)$ is biholomorphically equivalent to one of the pairs
$$
\begin{array}{cll}
\hbox{\rm (i)}& (\CC^N,\SU_N\ltimes \CC^N) & \hbox{{\rm (}here $\dim G=N^2+2N-1${\rm )}},\\
\vspace{-0.3cm}\\
\hbox{\rm (ii)} & (\CC^4,e^{\hbox{\tiny $i\RR$}}\Sp_2\ltimes \CC^4) & \hbox{{\rm (}here $N=4$ and $\dim G=N^2+3=19${\rm )}},\\
\vspace{-0.3cm}\\
\hbox{\rm (iii)} & (\CC^4,\Sp_2\ltimes \CC^4) & \hbox{{\rm (}here $N=4$ and $\dim G=N^2+2=18${\rm )}},\\
\vspace{-0.3cm}\\
\hbox{\rm (iv)} & (X_1\times X_2, G_1\times G_2) & \hbox{{\rm (}here $\dim G=N^2+2${\rm )}},
\end{array}
$$
where in {\rm (iv)} the manifold $X_1$ is one of $\CC^{N-1}$, $\CC\PP^{N-1}$, $\BB^{N-1}$, the manifold $X_2$ is one of $\CC$, $\CC\PP^1$, $\BB^1$, with the group $G_1$ being one of $\U_{N-1}\ltimes \CC^{N-1}$, $\PSU_N$, $\PSU_{N-1,1}$ and the group $G_2$ being one of $\U_1\ltimes \CC$, $\PSU_2$, $\PSU_{1,1}$, respectively.
\end{theorem}

\begin{remark}\label{overlapsrealcase} As we noted above (see (\ref{inequalityoverlap})), for $N=2,3,4$ the group dimension range that arises in the complex case overlaps with the group dimension range studied in the smooth situation in Sect.~\ref{riemanniancase}. Specifically, for $N=2$ the overlaps arise from \cite{Ish}, Theorem 4, cases $r=6,7,8$, as well as Theorem $\hbox{A}'$; for $N=3$ from \cite{KNa}, cases (4):(iv) and (5), p.~196; for $N=4$ from \cite{KNa}, cases (6):(iii) and (7):(iii), p.~197. It is not hard to observe that Theorem \ref{maxdimcompl} for $N=2$ and Theorem \ref{firstlacunacomplex} incorporate all these possibilities. At the same time, the proof of the theorems is independent of the approaches of \cite{Ish} and \cite{KNa}. Indeed, one can obtain Theorem \ref{firstlacunacomplex} either by observing that $X$ endowed with a $G$-invariant Hermitian metric is a Hermitian symmetric space (see \cite{Isa8}) or by a direct argument, which utilizes, in particular, the results of \cite{GK} and \cite{BDK} mentioned in Remark \ref{otherproofs} (see \cite{Isa5}).
\end{remark} 

As a consequence of Theorem \ref{firstlacunacomplex} we obtain:

\begin{corollary}\label{firstlacunacomplexcor} For $N\ge 3$, $N\ne 4$ a group $G$ satisfying
$$
N^2+3\le\dim G<N^2+2N-1
$$
cannot act properly on $X$. Furthermore, for $N=4$ there are no proper actions of groups of dimensions $20, 21, 22$.  
\end{corollary}

\noindent Corollary \ref{firstlacunacomplexcor} is analogous to Theorem \ref{firstlacuna} and shows that for $N\ge 3$ the set of the dimensions of the groups that can act properly on $X$ has a lacuna of size linear in $N$ located immediately below the second largest dimension $N^2+2N-1$.

We will now consider the case $\dim G=N^2+1$, with $N\ge 2$. In this situation, to describe the linear isotropy subgroups, one needs to determine connected closed subgroups of $\U_N$ of dimension $(N-1)^2$. They are given by the following lemma, which extends Lemma 2.1 of \cite{IKru1} (see \cite{Isa8}, \cite{Isa9}, \cite{IKru3}):

\begin{lemma}\label{un} \it Let $H\subset \U_N$ be a connected closed subgroup, with $N\ge 2$. Suppose that $\dim H=(N-1)^2$. Then $H$ is conjugate in $\U_N$ either to
\vspace{0.3cm}

\noindent {\rm (A)} $e^{\hbox{\tiny $i\RR$}}\SO_3(\RR)$ {\rm (}here $N=3${\rm )}; 
\vspace{0.3cm}

\noindent or to
\vspace{0.3cm}

\noindent {\rm (B)} $\SU_{N-1}\times \U_1$ realized by matrices in block-diagonal form, $N\ge 3$;
\vspace{0.3cm}

\noindent or, for some integers $k_1,k_2$, with $(k_1,k_2)=1, k_1>0$,  to
\vspace{0.3cm}

\noindent {\rm (C)} 
$$
H_{k_1,k_2}^N:=\left\{\left(\begin{array}{cc}
A & 0\\
0 & a
\end{array}\right)\in\U_N:A\in \U_{N-1},\, a\in
(\det A)^{\frac{k_2}{k_1}}\right\}.\label{mat}
$$
\end{lemma}

Now, we say that a pair $(X,G)$ (or the action of $G$ on $X$) is of type (A), (B) or (C), if for some $x\in X$ the group $LG_x^{\circ}$, when regarded as a subgroup of $\U_N$, is conjugate in $\U_N$ to the subgroup shown in (A), (B) or one of the subgroups defined in (C), respectively. Since the action of $G$ on $X$ is transitive, this definition is independent of the choice of $x\in X$. 

Actions of types (A) and (B) are not very hard to describe (although the latter are more difficult to deal with than the former). They are given by the following theorems (see \cite{Isa8}, \cite{Isa9}, \cite{IKru3}):

\begin{theorem}\label{classtype1} Any pair $(X,G)$ of type {\rm (A)} is biholomorphically equivalent to one of the pairs:
\vspace{0.3cm}

\noindent {\rm (i)} $(\mathscr{S},\SO(3,2)^{\circ})$, where 
$$
\mathscr{S}:=\Bigl\{(z_1,z_2,z_3)\in\CC^3: (\Im z_1)^2+(\Im z_2)^2-(\Im z_3)^2<0,\,\Im z_3>0\Bigr\}\footnote{Note that $\mathscr{S}$ is a tube representation of the symmetric classical domain of type $({\rm IV}_3)\simeq({\rm III}_2)$.}
$$
and $\Bihol(\mathscr{S})$ is identified with $\SO(3,2)^{\circ}$;
\vspace{0.3cm}

\noindent {\rm (ii)} $({\mathcal Q}_3,\SO_5(\RR))$, where ${\mathcal Q}_3$ is the projective quadric 
$$
{\mathcal Q}_3:=\left\{(z_1:z_2:z_3:z_4:z_5)\in\CC\PP^4: z_1^2+z_2^2+z_3^2+z_4^2+z_5^2=0\right\}
$$
and the compact real form of the complex group $\Bihol({\mathcal Q}_3)\simeq\SO_5(\CC)$ is identified with $\SO_5(\RR)$;
\vspace{0.3cm}

\noindent {\rm (iii)} $(\CC^3,e^{\hbox{\tiny $i\RR$}}\SO_3(\RR)\ltimes\CC^3)$.
\end{theorem}

\begin{theorem}\label{classtype2} Any pair $(X,G)$ of type {\rm (B)} is biholomorphically equivalent to a pair of the form $(\CC^{N-1}\times X',(\SU_{N-1}\ltimes\CC^{N-1})\times G')$, where the manifold $X'$ is one of $\CC$, $\CC\PP^1$, $\BB^1$, with the group $G'$ being one of $\U_1\ltimes\CC$, $\PSU_2$, $\PSU_{1,1}$, respectively.
\end{theorem}

We will now turn to actions of type (C). Describing such actions is the most challenging part of our classification, and we start with a large number of examples. Some of the examples can be naturally combined together and some of them form parametric families.
\vspace{0.5cm}

\noindent {\bf (i)} Here both the manifolds and the groups are represented as direct products.
\vspace{0.3cm}

{\bf (ia)} $X=X'\times\CC$, where $X'$ is one of $\CC^{N-1}$, $\CC\PP^{N-1}$, $\BB^{N-1}$ and\linebreak $G=G'\times\CC$, where $G'$ is one of the groups $\U_{N-1}\ltimes\CC^{N-1}$, $\PSU_{N}$, $\PSU_{N-1,1}$, respectively.
\vspace{0.3cm}

{\bf (ib)} $X=X'\times\CC^*$ and $G=G'\times \CC^*$, where $X'$ and $G'$ are as in\linebreak part {\bf (ia)}. 
\vspace{0.3cm}

{\bf (ic)} $X=X'\times\TT$ and $G=G'\times\Bihol(\TT)^{\circ}$, where $\TT$ is an elliptic curve and $X'$ and $G'$ are as in part {\bf (ia)}.
\vspace{0.3cm}

{\bf (id)} $X=X'\times{\mathcal P}_>$ and $G=G'\times (\RR\ltimes \RR)$, where ${\mathcal P}_>$ is the right half-plane $\left\{\xi\in\CC:\Re\xi>0\right\}$, the manifold $X'$ and the group $G'$ are as in part {\bf (ia)}, and $\RR\ltimes \RR$ is the group of all transformations of ${\mathcal P}_>$ of the form
$$
\xi\mapsto \lambda\xi+ia,\label{gofp}
$$
with $\lambda>0$, $a\in\RR$.
\vspace{0.5cm}

\noindent{\bf (ii)} Parts {\bf (iib)} and {\bf (iic)} below are obtained by passing to quotients in part {\bf (iia)}.

{\bf (iia)} $X=\BB^{N-1}\times\CC$ and $G$ consists of all maps of the form
$$
\begin{array}{lll}
z'&\mapsto&\displaystyle\frac{Az'+b}{cz'+d},\\
\vspace{0mm}&&\\
z_N&\mapsto&\displaystyle z_N+\ln(cz'+d)+a,
\end{array}
$$
where
$$
\left(\hspace{-0.1cm}\begin{array}{cc}
A& b\\
c& d
\end{array}
\hspace{-0.1cm}\right)
\in \SU_{N-1,1},\quad z':=\left(\hspace{-0.1cm}\begin{array}{c}z_1\\\vdots\\z_{N-1}\end{array}\hspace{-0.1cm}\right)\label{autballn-1}
$$
and $a\in\CC$. In fact, for $t\in\CC$ one can consider the following family of groups acting on $X$:
\begin{equation}
\begin{array}{lll}
z'&\mapsto&\displaystyle\frac{Az'+b}{cz'+d},\\
\vspace{0mm}&&\\
z_N&\mapsto&\displaystyle z_N+t\ln(cz'+d)+a,
\end{array}\label{groupbcs}
\end{equation}
where $A,a,b,c,d$ are as above. Example {\bf (ia)} for $X'=\BB^{N-1}$ is included in this family for $t=0$. If $t\ne 0$, conjugating group (\ref{groupbcs}) in $\Bihol(X)$ by the automorphism
\begin{equation}
\begin{array}{lll}
z'&\mapsto & z',\\
z_N&\mapsto & z_N/t,
\end{array}\label{divs}
\end{equation}
we can assume that $t=1$.  
\vspace{0.3cm}

{\bf (iib)} $X=\BB^{N-1}\times\CC^*$ and for a fixed $t\in\CC^*$ the group $G$ consists of all maps of the form
\begin{equation}
\begin{array}{lll}
z'&\mapsto&\displaystyle\frac{Az'+b}{cz'+d},\\
\vspace{0mm}&&\\
z_N&\mapsto&\displaystyle\chi (cz'+d)^tz_N,
\end{array}\label{groupbcstar}
\end{equation}
where $A,b,c,d$ are as in {\bf (iia)} and $\chi\in\CC^*$. Example {\bf (ib)} for $X'=\BB^{N-1}$ can be included in this family for $t=0$. Clearly, (\ref{groupbcstar}) is obtained from (\ref{groupbcs}) by passing to a quotient in $z_N$.
\vspace{0.3cm}

{\bf (iic)} $X=\BB^{N-1}\times\TT$, where $\TT$ is an elliptic curve, and for a fixed $t\in\CC^*$ the group $G$ consists of all maps of the form
\begin{equation}
\begin{array}{lll}
z'&\mapsto&\displaystyle\frac{Az'+b}{cz'+d},\\
\vspace{0mm}&&\\

[z_N]&\mapsto&\displaystyle
\left[\chi (cz'+d)^tz_N\right].
\end{array}\label{groupbtorus}
\end{equation}
Here $A,b,c,d,\chi$ are as in {\bf (iib)}, $\TT$ is obtained from $\CC^*$ by taking the quotient with respect to the equivalence relation $z_N\sim hz_N$, for some $h\in\CC^*$, $|h|\ne 1$, and $[z_N]\in\TT$ is the equivalence class of a point $z_N\in\CC^*$. Example {\bf (ic)} for $X'=\BB^{N-1}$ can be included in this family for $t=0$. Clearly, (\ref{groupbtorus}) arises from (\ref{groupbcstar}) by passing to the quotient in $z_N$ specified above.
\vspace{0.5cm}

\noindent{\bf (iii)} Part {\bf (iiib)} below is obtained by passing to a quotient in\linebreak part {\bf (iiia)}.  
\vspace{0.3cm}

{\bf (iiia)} $X=\CC^N$ and $G$ consists of all maps of the form
$$
\begin{array}{lll}
z'&\mapsto & e^{\hbox{\tiny Re}\,b}Uz'+a,\\
z_N&\mapsto & z_N+b,
\end{array}
$$
where $U\in \U_{N-1}$, $a\in\CC^{N-1}$, $b\in\CC$. In fact, for $t\in\CC$ one can consider the following family of groups acting on $X$:
\begin{equation}
\begin{array}{lll}
z'&\mapsto & e^{\hbox{\tiny Re}\,(t b)}Uz'+a,\\
z_N&\mapsto & z_N+b,
\end{array}\label{g3cns}
\end{equation}
where $U$, $a$, $b$ are as above. Example {\bf (ia)} for $X'=\CC^{N-1}$ is included in this family for $t=0$. If $t\ne 0$, conjugating group (\ref{g3cns}) in $\Bihol(X)$ by the automorphism
$$
\begin{array}{lll}
z'&\mapsto & z',\\
z_N&\mapsto & t z_N,
\end{array}
$$
we can assume that $t=1$.  
\vspace{0.3cm}

{\bf (iiib)} $X=\CC^{N-1}\times\CC^*$ and for a fixed $t\in\RR^*$ the group $G$ consists of all maps of the form
\begin{equation}
\begin{array}{lll}
z'&\mapsto & e^{t\,\hbox{\tiny Re}\,b}Uz'+a,\\
z_N&\mapsto & e^{b}z_N,
\end{array}\label{iiib}
\end{equation}
where $U, a, b$ are as in {\bf (iiia)}. Example {\bf (ib)} for $X'=\CC^{N-1}$ can be included in this family for $t=0$. Clearly, (\ref{iiib}) is obtained from (\ref{g3cns}) for $t\in\RR^*$ by passing to a quotient in $z_N$. 
\vspace{0.5cm}

\noindent{\bf (iv)} Parts {\bf (ivb)} and {\bf (ivc)} below are obtained by passing to quotients in part {\bf (iva)}.
\vspace{0.3cm}

{\bf (iva)} $X=\CC^N$ and $G$ consists of all maps of the form
$$
\begin{array}{lll}
z' & \mapsto & Uz'+a,\\
z_N & \mapsto & z_N+\langle Uz',a\rangle+b,
\end{array}
$$
where $U\in \U_{N-1}$, $a\in\CC^{N-1}$, $b\in\CC$, and $\langle\cdot\,,\cdot\rangle$ is the Hermitian scalar product on $\CC^{N-1}$. In fact, for $t\in\CC$ one can consider the following family of groups acting on $X$:
\begin{equation}
\begin{array}{lll}
z' & \mapsto & Uz'+a,\\
z_N & \mapsto & z_N+t\langle Uz',a\rangle+b,
\end{array}\label{g4cns}
\end{equation}
where $U$, $a$, $b$ are as above. Example {\bf (ia)} for $X'=\CC^{N-1}$ is included in this family for $t=0$. If $t\ne 0$, conjugating group (\ref{g4cns}) in $\Bihol(X)$ by automorphism (\ref{divs}), we can assume that $t=1$.  
\vspace{0.3cm}

{\bf (ivb)} $X=\CC^{N-1}\times\CC^*$ and for a fixed $0\le\tau<2\pi$ the group $G$ consists of all maps of the form
\begin{equation}
\begin{array}{lll}
z' & \mapsto & Uz'+a,\\
z_N & \mapsto &\chi\exp\Bigl(e^{i\tau}\langle Uz',a\rangle\Bigr)z_N,
\end{array}\label{groupcstartau}
\end{equation}
where $U,a$ are as in {\bf (iva)} and $\chi\in\CC^*$. In fact, for $t\in\CC$ one can consider the following family of groups acting on $X$:
\begin{equation}
\begin{array}{lll}
z' & \mapsto & Uz'+a,\\
z_N & \mapsto &\chi\exp\Bigl(t\langle Uz',a\rangle\Bigr)z_N,
\end{array}\label{g4cnsstar}
\end{equation}
where $U,a,\chi$ are as above. Example {\bf (ib)} for $X'=\CC^{N-1}$ is included in this family for $t=0$. Formula (\ref{g4cnsstar}) is obtained from (\ref{g4cns}) by passing to a quotient in $z_N$. Furthermore, if $t\ne 0$, conjugating group (\ref{g4cnsstar}) in $\Bihol(X)$ by the automorphism
$$
\begin{array}{lll}
z'&\mapsto & \sqrt{|t|}z',\\
z_N&\mapsto & z_N,
\end{array}
$$
we obtain the group defined in (\ref{groupcstartau}) for $\tau=\hbox{arg}\,t$.
\vspace{0.3cm}

{\bf (ivc)} $X=\CC^{N-1}\times\TT$, where $\TT$ is an elliptic curve  obtained from $\CC^*$ as in {\bf (iic)}, and for a fixed $0\le\tau<2\pi$ the group $G$ consists of all maps of the form
\begin{equation}
\begin{array}{lll}
z' & \mapsto & Uz'+a,\\

[z_N] & \mapsto &\left[\chi\exp\Bigl(e^{i\tau}\langle Uz',a\rangle\Bigr)z_N\right],
\end{array}\label{grouptorustau}
\end{equation}
with $U,a,\chi$ being as in {\bf (ivb)}. In fact, for $t\in\CC$ one can consider the following family of groups acting on $X$:
\begin{equation}
\begin{array}{lll}
z' & \mapsto & Uz'+a,\\

[z_N] & \mapsto &\left[\chi\exp\Bigl(t\langle Uz',a\rangle\Bigr)z_N\right],
\end{array}\label{g4cnst}
\end{equation}
where $U,a,\chi$ are as above. Example {\bf (ic)} for $X'=\CC^{N-1}$ is included in this family for $t=0$. Clearly, formula (\ref{g4cnst}) is derived from (\ref{g4cnsstar}) by passing to a quotient in the variable $z_N$. Furthermore, if $t\ne 0$, conjugating group (\ref{g4cnst}) in $\Bihol(X)$ by the automorphism
$$
\begin{array}{lll}
z'&\mapsto & \sqrt{|t|}z',\\
\xi &\mapsto & \xi,
\end{array}
$$
where $\xi\in\TT$, we obtain the group defined in (\ref{grouptorustau}) for $\tau=\hbox{arg}\,t$.
\vspace{0.5cm}

\noindent{\bf (v)} $X=\CC^{N-1}\times{\mathcal P}_>$ and for a fixed $t\in\RR^*$ the group $G$ consists of all maps of the form
$$
\begin{array}{lll}
z'&\mapsto & \lambda^tUz'+a,\\
z_N&\mapsto & \lambda z_N+ib,
\end{array}\label{groupv}
$$
where $U\in \U_{N-1}$, $a\in\CC^{N-1}$, $b\in\RR$, $\lambda>0$. Example {\bf (id)} for $X'=\CC^{N-1}$ can be included in this family for $t=0$.
\vspace{0.5cm}

\noindent {\bf (vi)} $X=\CC^N$ and for fixed $k_1,k_2\in\ZZ$, $(k_1,k_2)=1$, $k_1>0$, $k_2\ne 0$ the group $G$ consists of all maps of the form 
$$
z\mapsto Uz+a,
$$
where $a\in\CC^N$ and $U\in H_{k_1,k_2}^N$ (see (C) in Lemma \ref{un}). Example {\bf (ia)} for $X'=\CC^{N-1}$ can be included in this family for $k_2=0$.  
\vspace{0.5cm}

\noindent{\bf (vii)} Part {\bf (viib)} below is obtained by passing to a quotient in\linebreak part {\bf (viia)}.  
\vspace{0.3cm}

{\bf (viia)} $X=\CC^{N{}*{}}/\ZZ_l$, with $\CC^{N{}*{}}:=\CC^N\setminus\{0\}$, $l\in\NN$, and $G$ consists of all elements of $\Bihol(X)$ induced by transformations of $\CC^{N{}*{}}$ of the form
\begin{equation}
z\mapsto\lambda Uz,\label{gviiia}
\end{equation}
where $U\in \U_N$ and $\lambda>0$.
\vspace{0.3cm}

{\bf (viib)} $X=M_h/\ZZ_l$, with $M_h$ being the Hopf manifold $\CC^{N{}*{}}/\{z{}\sim{}hz\}$, for $h\in\CC^*$, $|h|\ne 1$, and $G$ consists of all elements of $\Bihol(X)$ induced by transformations of $\CC^{N{}*{}}$ of the form (\ref{gviiia}). 
\vspace{0.5cm}

\noindent{\bf (viii)} Let $G$ be the group of all affine maps of $\CC^N$ of the form
$$
\begin{array}{lll}
z' & \mapsto & \lambda Uz'+a,\\
z_N & \mapsto & \lambda^2z_N+2\lambda\langle Uz',a\rangle+|a|^2+ib,
\end{array}
$$
where $U\in \U_{N-1}$, $a\in\CC^{N-1}$, $b\in\RR$, $\lambda>0$. Here the manifolds are the open orbits of the action of $G$ on $\CC^N$.
\vspace{0.3cm}

{\bf (viiia)}
$$
X=\left\{z\in\CC^N: \Re z_N>|z'|^2\right\}.
$$
Observe that this domain is biholomorphic to $\BB^N$. 
\vspace{0.3cm}

{\bf (viiib)}
$$
X=\left\{z\in\CC^N: \Re z_N<|z'|^2\right\}.
$$
Observe that this domain is biholomorphic to $\CC\PP^N\setminus(\overline{\BB^N}\cup L)$, where $L$ is a complex hyperplane tangent to $\partial\BB^N$.
\vspace{0.5cm}

\noindent {\bf (ix)} Here $N=2$, $X=\BB^1\times\CC$ and $G$ consists of all maps of the form
$$
\begin{array}{lll}
z_1 & \mapsto & \displaystyle\frac{az_1+b}{\overline{b}z_1+\overline{a}},\\
\vspace{0.1cm}\\
z_2 & \mapsto & \displaystyle\frac{z_2+cz_1+\overline{c}}{\overline{b}z_1+\overline{a}},\\
\end{array}
$$
where $a,b\in\CC$, $|a|^2-|b|^2=1$, $c\in\CC$.
\vspace{0.5cm}

\noindent{\bf (x)} Here $N=3$, $X=\CC\PP^3$ and $G$ consists of all maps of the form 
$$
\zeta\mapsto U\zeta,
$$
where $U\in \Sp_2$ and $\zeta$ is a point in $\CC\PP^3$ written as a column vector of four homogeneous coordinates.
\vspace{0.5cm}

\noindent{\bf (xi)} Let $N=3$ and write a point $\zeta\in\CC\PP^3$ in homogeneous coordinates as $\zeta=(\zeta':\zeta'')$, with $\zeta'=(\zeta'_1:\zeta'_2)$, $\zeta''=(\zeta''_1:\zeta''_2)$. Now, set\linebreak $X=\CC\PP^3\setminus\{\zeta''=0\}$ and let $G$ be the group of all maps of the form
$$
\begin{array}{lll}
\zeta'&\mapsto&U\zeta'+A\zeta'',\\
\zeta''&\mapsto& V\zeta'',
\end{array}\label{gcp3cp1}
$$
where $U,V\in \SU_2$,   
$$
A=\left(
\begin{array}{lr}
a & i\overline{b}\\
b & -i\overline{a}
\end{array}
\right),
$$
for some $a,b\in\CC$, and $\zeta'$, $\zeta''$ are written as column vectors. 
\vspace{0.5cm}

\noindent{\bf (xii)} Here $N=3$, $X=\CC^3$ and $G$ consists of all maps of the form
$$
\begin{array}{lll}
z'&\mapsto&Uz'+a,\\
z_3&\mapsto&\hbox{det}\,U\,z_3+\left[
\left(
\begin{array}{ll}
0 & 1\\
-1 & 0
\end{array}
\right)
Uz'\right]\cdot a+b,
\end{array}
$$
where $U\in \U_2$, $a\in\CC^2$, $b\in\CC$, and $\cdot$ is the bilinear form on $\CC^2$ defined as follows: $(u_1,u_2)\cdot(v_1,v_2):=u_1v_1+u_2v_2$ for all $(u_1,u_2), (v_1,v_2)\in\CC^2$. 
\vspace{0.3cm}

We can now state the main result of this section (see \cite{IKru3}):

\begin{theorem}\label{classtype3} Any pair $(X,G)$ of type {\rm (C)} is biholomorphically equivalent to one of the pairs listed in {\bf (i)--(xii)}.
\end{theorem}

Studying proper actions becomes harder as the group dimension decreases, and there exists no complete classification of pairs $(X,G)$ for any prescribed $\dim G$ not exceeding $N^2$. There are only partial results in the cases $\dim G=N^2$, $\dim G=N^2-1$, which we will now briefly mention. 

Firstly, for $G$ isomorphic to either $\U_N$ (here $\dim G=N^2$) or $\SU_N$ (here $\dim G=N^2-1$), with $N\ge 2$, classifications of pairs $(X,G)$ were obtained in \cite{IKru1}, \cite{IKru2} (see also \cite{Isa6}, Sects. 6.2, 6.3). Secondly, suppose that $\dim G=N^2$. As the action of $G$ may not be transitive, one has to consider the transitive and non-transitive cases separately. To understand transitive actions, one needs to have a classification of closed connected subgroups of $\U_N$ of dimension $N^2-2N$, with $N\ge 2$, and such subgroups were to some extent determined in Lemma 4.2 of \cite{IKru1}. On the other hand, for non-transitive actions we have the following description of the orbits (see \cite{Isa6}, Proposition 6.3):

\begin{proposition}\label{main1prop} Let $N\ge 2$, $\dim G=N^2$ and the action of $G$ on $X$ be non-transitive. Fix $x\in X$ and set $V:=T_x(Gx)$. Then
\vspace{-0.3cm}\\

\noindent {\rm (i)} either the orbit $Gx$ is a real or complex closed hypersurface in $X$, or $x$ is a fixed point of the action;
\vspace{-0.3cm}\\

\noindent {\rm (ii)} if $Gx$ is a real hypersurface, it is either spherical or Levi-flat, and in the latter case it is foliated by complex manifolds each of which is biholomorphic to one of $\CC^{N-1}$, $\CC\PP^{N-1}$, $\BB^{N-1}$; there exist coordinates in $T_x(X)$ such that, with respect to the orthogonal decomposition
$$
T_x(X)=(V\cap iV)\oplus(V\cap iV)^{\bot},
$$
the group $LG_x$ is $\U_{N-1}\times K$, where $K$ is either $\ZZ_2$ or trivial, and the former can only occur if $Gx$ is Levi-flat; if $Gx$ is Levi-flat and ${\mathfrak F}(x)$ is the leaf of the foliation passing through $x$, then a biholomorphism between ${\mathfrak F}(x)$ and one of $\CC^{N-1}$, $\CC\PP^{N-1}$, $\BB^{N-1}$ can be chosen so that it transforms $G_x/\alpha_x^{-1}(K)$, regarded as a subgroup of $\Bihol({\mathfrak F}(x))$, into $\U_{N-1}$ considered as a subgroup of one of $\U_{N-1}\ltimes\CC^{N-1}$, $\PSU_N$, $\PSU_{N-1,1}$, respectively;
\vspace{-0.3cm}\\

\noindent {\rm (iii)} if $Gx$ is a complex hypersurface, it is biholomorphic to one of $\CC^{N-1}$, $\CC\PP^{N-1}$, $\BB^{N-1}$; there exist coordinates in $T_x(X)$ such that, with respect to the orthogonal decomposition $T_x(X)=V\oplus V^{\bot}$, we have
\begin{equation}
LG_x=\U_{N-1}\times\U_1,\label{linisotrproduct}
\end{equation}
and a biholomorphism between $Gx$ and one of $\CC^{N-1}$, $\CC\PP^{N-1}$, $\BB^{N-1}$ can be chosen so that it transforms $G/\alpha_{x}^{-1}(\U_1)$ {\rm(}where $\U_1$ is the second factor in {\rm(\ref{linisotrproduct})}{\rm)}, regarded as a subgroup of $\Bihol(Gx)$, into one of $\U_{N-1}\ltimes\CC^{N-1}$, $\PSU_N$, $\PSU_{N-1,1}$ respectively; there are at most two complex hypersurface orbits in $X$;
\vspace{-0.3cm}\\

\noindent {\rm (iv)} if $x$ is a fixed point of the action, then $X$ is biholomorphic to one of $\CC^N$, $\CC\PP^N$, $\BB^N$, and the biholomorphism can be chosen so that it transforms $G$ into $\U_N$ regarded as a subgroup of one of $\U_N\ltimes\CC^N$, $\PSU_{N+1}$, $\PSU_{N,1}$, respectively.

\end{proposition}

By Proposition \ref{main1prop}, in order to obtain a classification for group dimension $N^2$ in the non-transitive case, one needs to \lq\lq glue\rq\rq\hspace{-0.01cm} together copies of $\CC^{N-1}$, $\CC\PP^{N-1}$, $\BB^{N-1}$, spherical hypersurfaces, as well as Levi-flat hypersurfaces foliated by copies of $\CC^{N-1}$, $\CC\PP^{N-1}$, $\BB^{N-1}$. This task appears to be rather hard in general. However, as we will see in the next section, it can be successfully completed for Kobayashi-hyperbolic manifolds. Specifically, there exists a biholomorphic classification of all Kobayashi-hyperbolic manifolds $X$ in the following two situations: $\dim\Bihol(X)=N^2$ and $\dim\Bihol(X)=N^2-1$ (below we will observe that the action of the group $\Bihol(X)$ is proper). Kobayashi-hyperbolic manifolds form an important class, and investigating manifolds with prescribed dimension of the automorphism group in this class was in fact one of the main motivations for our study of general proper actions. 

\section{Kobayashi-hyperbolic manifolds}\label{kobayashihyperboliccase}
\setcounter{equation}{0}

In this section, $X$ continues to denote a connected complex paracompact manifold of complex dimension $N$. We start by introducing a certain pseudodistance on $X$, called the Kobayashi pseudodistance, and stating some of its basic properties. For more detail the reader is referred to \cite{Ko2}, Chap.~3.

First, recall that the {\it Poincar\'e distance}\, on $\BB^1$ is given by
$$
\rho(z,w):=\displaystyle\frac{1}{2}\ln\frac{1+\left|\displaystyle\frac{z-w}{1-\overline{z}w}\right|}{1-\left|\displaystyle\frac{z-w}{1-\overline{z}w}\right|},\quad z,w\in\BB^1.
$$
The {\it Kobayashi pseudodistance}\, on $X$ is then defined as follows:
$$
K_X(x,y):=\inf\sum_{j=1}^m\rho(z_j,w_j),\quad x,y\in X,
$$
where the $\inf$ is taken over all $m\in\NN$, all pairs of points $\{z_j,w_j\}_{j=1,\dots,m}$ in $\BB^1$ and all collections of holomorphic maps $\{f_j\}_{j=1,\dots,m}$ from $\BB^1$ into $X$ such that $f_1(z_1)=x$, $f_m(w_m)=y$, and $f_j(w_j)=f_{j+1}(z_{j+1})$ for\linebreak $j=1,\dots,m-1$. Clearly, one has $K_{\BB^1}=\rho$. 

It is straightforward to verify that $K_X$ is indeed a pseudodistance on $X$. Furthermore, it does not increase under holomorphic maps, i.e., for any holomorphic map $f$ between two manifolds $X_1$ and $X_2$ we have
\begin{equation}
K_{X_{{}_2}}\left(f(x),f(y)\right)\le K_{X_{{}_1}}(x,y),\quad x,y\in X_1.\label{kobayashinonincrease}
\end{equation}
In particular, $K_X$ is $\Bihol(X)$-invariant.

One can easily produce examples where $K_X$ is not a distance, that is, $K_X(x,y)=0$ for some $x,y\in X$, $x\ne y$. For instance, the definition of $K_X$ immediately implies $K_{\CC}\equiv 0$. Hence, for any manifold $X$ admitting an {\it entire curve}, i.e., a non-constant holomorphic map $f:\CC\to X$, the pseudodistance $K_X$ is not a distance as by (\ref{kobayashinonincrease}) one has $K_X\equiv 0$ on $f(\CC)$. A manifold $X$ for which $K_X$ is a distance is called {\it Kobayashi-hyperbolic}. Note that a Riemann surface is Kobayashi-hyperbolic if and only if it is {\it hyperbolic}, i.e., its universal cover is biholomorphic\linebreak to $\BB^1$.

It is easy to see that for any manifold $X$ the function $K_X$ is continuous on $X\times X$, hence Theorem \ref{complexmetrics} implies that if $X$ is Kobayashi-hyperbolic, the action of the group $\Bihol(X)$ on $X$ is proper. We also note in passing that if $X$ is Kobayashi-hyperbolic, $K_X$ induces the topology of $X$ (see \cite{Ro}). This fact, as well as the continuity of $K_X$ for any $X$, holds true for complex spaces as well (see \cite{Ba1} and \cite{Ba2}).

The class of Kobayashi-hyperbolic manifolds is quite large. In fact, it cannot be expanded by means of considering other pseudodistances that do not increase under holomorphic maps. Indeed, the Kobayashi pseudodistance possesses a maximality property, which we will now formulate. 

First, an assignment of a pseudodistance ${\mathtt d}_X$ to every manifold $X$ is called a {\it Schwarz-Pick system}\, if
\vspace{-0.3cm}\\

\noindent (i) ${\mathtt d}_{\BB^1}=\rho$;
\vspace{0.3cm}

\noindent (ii) for any holomorphic map $f$ between any two manifolds $X_1$, $X_2$ one has
$$
{\mathtt d}_{X_{{}_2}}\left(f(x),f(y)\right)\le {\mathtt d}_{X_{{}_1}}(x,y),\quad x,y\in X_1.
$$
Now, the maximality property of the Kobayashi pseudodistance is the fact that $\{K_X\}$ is the largest Schwarz-Pick system, that is, one has $K_X\ge {\mathtt d}_X$ for every $X$ and any Schwarz-Pick system $\{{\mathtt d}_X\}$ (see, e.g., \cite{Din}, p.~52). In particular, if for some Schwarz-Pick system $\{{\mathtt d}_X\}$ and some manifold $X_0$ the pseudodistance ${\mathtt d}_{X_{{}_0}}$ is in fact a distance, then $X_0$ is Kobayashi-hyperbolic. 

Using the above maximality property, one can compare $K_X$ with the classical {\it Carath\'eodory pseudodistance}\, defined as follows:
$$
C_X(x,y):=\sup_f\rho(f(x),f(y)),\quad x,y\in X,
$$
where the $\sup$ is taken over all holomorphic maps $f:X\ra\BB^1$. It is easy to show that $\{C_X\}$ is a Schwarz-Pick system; in fact, it is the smallest one in the sense that $C_{X}\le {\mathtt d}_{X}$ for every manifold $X$ and any Schwarz-Pick system $\{{\mathtt d}_X\}$ (see, e.g., \cite{Din}, Proposition 4.2). Clearly, $C_X$ is a distance on any manifold $X$ whose points are separated by bounded holomorphic functions (for instance, on bounded domains in complex space). The maximality property of the Kobayashi pseudodistance (or, alternatively, the minimality property of the Carath\'eodory pseudodistance) then yields that all such manifolds are Kobayashi-hyperbolic. There are, however, many more Kobayashi-hyperbolic manifolds than those on which the Carath\'eodory pseudodistance does not degenerate. Indeed, for any compact manifold $X$ one has $C_X\equiv 0$, but nevertheless, if $X$ does not admit entire curves, it is Kobayashi-hyperbolic (see \cite{Br}).

From now on, we let $X$ to be Kobayashi-hyperbolic and work with the connected group $G(X):=\Bihol(X)^{\circ}$. This section contains complete biholomorphic classifications of manifolds $X$ in the cases\linebreak $\dim G(X)=N^2$, $\dim G(X)=N^2-1$, which extend results of Sect.~\ref{complexcase} for $G=G(X)$. Everywhere below we assume $N\ge 2$. 

Before proceeding, we will explain why $N^2-1$ is the lowest automorphism group dimension for which one can hope to obtain an explicit classification. Indeed, let us focus for the moment on {\it Reinhardt domains}\, in $\CC^N$, i.e., the domains invariant under the torus action $z_j\mapsto e^{i\alpha_j}z_j$, $\alpha_j\in\RR$, $j=1,\dots,N$. \lq\lq Most\rq\rq\hspace{-0.01cm} Reinhardt domains have no biholomorphic automorphisms other than those induced by this action and hence have an $N$-dimensional automorphism group. In particular, if $D$ is a \lq\lq generic\rq\rq\hspace{-0.01cm}\ Kobayashi-hyperbolic Reinhardt domain in $\CC^2$, then $\dim G(D)=2=N^2-2$. Such Reinhardt domains cannot be explicitly described up to biholomorphic equivalence.

To make this argument more precise, choose a subset $Q$ of the first quadrant $\RR^2_{+}:=\{(x,y)\in\RR^2: x\ge 0, y\ge 0\}$ so that the associated set in $\CC^2$
$$
D_Q:=\left\{(z_1,z_2)\in\CC^2: (|z_1|, |z_2|)\in Q\right\}
$$
is a smoothly bounded Reinhardt domain and contains the origin. Furthermore, it is easy to choose $Q$ satisfying the requirement that the group $\Bihol(Q_Q)$ be compact. Indeed, in \cite{FIK} all smoothly bounded Reinhardt domains with non-compact automorphism group were found, and all of them turned out to contain the origin. On the other hand, by \cite{Su}, two bounded Reinhardt domains containing the origin are biholomorphic to each other if and only if one is obtained from the other
by dilating and permuting coordinates. Therefore, one can easily make a choice of $Q$ so that $D_Q$ is not biholomorphic to any of the domains listed in \cite{FIK} thus ensuring that $\Bihol(D_Q)$ is compact. It then follows from the explicit description of the automorphism groups of bounded Reinhardt domains (see \cite{Kr} and \cite{Sh}) that $G(D_Q)$ is isomorphic to either $\U_2$ or $\U_1\times \U_1$. However, by \cite{GIK}, there does not exist a Kobayashi-hyperbolic Reinhardt domain in $\CC^2$ containing the origin for which the automorphism group is compact and 4-dimensional. Hence, $G(D_Q)$ is in fact isomorphic to $\U_1\times \U_1$ and thus $\dim G(D_Q)=2$. The freedom in choosing a set $Q\subset\RR^2_{+}$ that satisfies the above requirements is very substantial, and by varying $Q$ one can produce a family of pairwise biholomorphically non-equivalent Reinhardt domains $D_Q$ that cannot be described by any reasonable formulas.

The above argument shows that there is no explicit classification of Kobayashi-hyper\-bolic manifolds $X$ with $\dim G(X)=N^2-2$ if one insists on accommodating all dimensions $N\ge 2$. However, some acceptable classifications may exist for $N\ge 3$, $\dim G(X)=N^2-2$, as well as for particular values of $N$ and $\dim G(X)$, with $\dim G(X)< N^2-2$ (see, e.g., \cite{GIK} for a study of Reinhardt domains from the point of view of automorphism group dimensions).

We will now proceed with our analysis of the cases $\dim G(X)=N^2$, $\dim G(X)=N^2-1$. First, in the following theorem we consider the situation when $X$ is {\it homogeneous}, i.e., when the action of $G(X)$ is transitive (see \cite{Isa1}; \cite{Isa2}; \cite{Isa4}; \cite{Isa6}, Theorem 2.2):

\begin{theorem}\label{fullclasshom} Let $X$ be homogeneous. Then the following holds:
\vspace{-0.1cm}\\

\noindent{\rm (i)} if $\dim G(X)=N^2$, then either $N=3$ and $X$ is biholomorphic to $\BB^1\times\BB^1\times\BB^1$, or $N=4$ and $X$ is biholomorphic to $\BB^2\times \BB^2$;
\vspace{-0.1cm}\\

\noindent{\rm (ii)} if $\dim G(X)=N^2-1$, then $N=4$ and $X$ is biholomorphic to
$$
\begin{array}{ll}
\mathscr{T}:=\Bigl\{(z_1,z_2,z_3,z_4)\in\CC^4:&(\Im z_1)^2+(\Im z_2)^2+
\\&(\Im z_3)^2-(\Im z_4)^2<0,\,\Im z_4>0\Bigr\}.
\end{array}\hspace{-0.3cm}\footnote{Note that $\mathscr{T}$ is a tube realization of the symmetric classical domain of type $({\rm IV}_4)\simeq({\rm I}_{2,2})$ and $G(\mathscr{T})\simeq\SO(4,2)^{\circ}/\ZZ_2$.}
$$
\end{theorem}

While the above result includes only three domains, a complete description of general complex manifolds that admit transitive pro\-per actions of groups of dimension $N^2$ is probably much larger. Indeed, notice for comparison that the only Kobayashi-hyperbolic manifold with $(N^2+1)$-dimensional automorphism group on the very long list that appears in Theorems \ref{classtype1}, \ref{classtype2}, \ref{classtype3} is the tube domain $\mathscr{S}\subset\CC^3$. The only other Kobayashi-hyperbolic manifolds on the list are $\BB^{N-1}\times\BB^1$ and $\BB^N$, but one has $\dim G(\BB^{N-1}\times\BB^1)=N^2+2$ and $\dim G(\BB^N)=N^2+2N$.

We will now briefly explain the idea of the proof of Theorem \ref{fullclasshom}. Since $X$ is homogeneous, by \cite{N}, \cite{VGP-S} it is biholomorphic to a {\it Siegel domain of the second kind}, i.e., a domain of the form
$$
\left\{(z',z'')\in\CC^{N-k}\times\CC^k: \Im z''-F(z',z')\in C\right\}
$$
for an integer $1\le k\le N$, an open convex cone $C\subset\RR^k$ not containing an entire line, and a $\CC^k$-valued Hermitian form $F$ defined on $\CC^{N-k}\times\CC^{N-k}$ such that $F(z',z')\in\overline{C}\setminus\{0\}$ for all non-zero $z'\in\CC^{N-k}$. Every domain of this form is unbounded by definition, but in fact has a bounded realization and hence is Kobayashi-hyperbolic. It then remains to determine all homogeneous Siegel domains of the second kind, say $D$, with either $\dim G(D)=N^2$ or $\dim G(D)=N^2-1$. This can be achieved by utilizing an explicit description of the Lie algebra of the group $G(D)$ (see \cite{Sa}, pp.~211--219) and leads to the classification stated in the theorem.

Next, we consider the case when $\dim G(X)=N^2$ and $X$ is not homogeneous. For Kobayashi-hyperbolic manifolds, the statement of Proposition \ref{main1prop} significantly simplifies as one must exclude $\CC^{N-1}$, $\CC\PP^{N-1}$ from parts (ii), (iii), and part (iv) becomes irrelevant. Furthermore, one can show that the existence of two complex hypersurface orbits contradicts Kobayashi-hyperbolicity. Then a thorough analysis of the strongly pseudoconvex as well as Levi-flat orbits allows us to \lq\lq glue\rq\rq\hspace{-0.01cm} them together, possibly \lq\lq adding\rq\rq\hspace{-0.01cm} at most one complex hypersurface orbit, to obtain the following result (see \cite{Isa1}; \cite{Isa4}; \cite{Isa6},\linebreak Theorem 3.1):

\begin{theorem}\label{main} Let $X$ be non-homogeneous and $\dim G(X)=N^2$. Then $X$ is biholomorphic to one and only one of the manifolds listed below:
\vspace{-0.1cm}\\

\noindent{\rm (i)} $\{z\in\CC^N: r<||z||<1\}/\ZZ_l$, $0\le r<1$, $l\in\NN$;
\vspace{-0.1cm}\\

\noindent{\rm (ii)} $\Bigl\{(z',z_N)\in\CC^{N-1}\times\CC: ||z'||^2+|z_N|^{\theta}<1\Bigr\}$, $\theta>0$, $\theta\ne 2$;
\vspace{-0.1cm}\\

\noindent{\rm (iii)} $\Bigl\{(z',z_N)\in\CC^{N-1}\times\CC:||z'||<1,\, |z_N|<(1-||z'||^2)^{\theta}\Bigr\}$, $\theta<0$;
\vspace{-0.1cm}\\

\noindent{\rm (iv)} $\Bigl\{(z',z_N)\in\CC^{N-1}\times\CC: ||z'||<1,\, r(1-||z'||^2)^{\theta}<|z_N|<\linebreak(1-||z'||^2)^{\theta}\Bigr\}$, with either $\theta\ge 0$, $0\le r<1$, or $\theta<0$, $r=0$;
\vspace{-0.1cm}\\

\noindent{\rm (v)} $\Bigl\{(z',z_N)\in\CC^{N-1}\times\CC: r\exp\left({\theta||z'||^2}\right)<|z_N|<\exp\left({\theta||z'||^2}\right)\Bigr\}$, with either $\theta=1$, $0<r<1$, or $\theta=-1$, $r=0$;
\vspace{-0.1cm}\\

\noindent{\rm (vi)} $\Bigl\{(z',z_N)\in\CC^{N-1}\times\CC: ||z'||<1,\,r(1-||z'||^2)^{\theta}<\exp\left(\Re z_N\right)<(1-||z'||^2)^{\theta}\Bigr\}$, with either $\theta=1$, $0\le r<1$, or $\theta=-1$, $r=0$;
\vspace{-0.1cm}\\

\noindent{\rm (vii)} $\Bigl\{(z',z_N)\in\CC^{N-1}\times\CC: -1+||z'||^2<\Re z_N<||z'||^2\Bigr\}$.
\end{theorem}

\noindent We note that for each manifold provided by Theorem \ref{main}, determining the $N^2$-dimensional connected identity component of its automorphism group explicitly is not hard (see, e.g., \cite{Isa6}, p.~31).

\begin{remark}\label{hostorical} For the first time, the situation $\dim G(X)=N^2$ was considered in article \cite{GIK}, where we found a classification in the case of Reinhardt domains in $\CC^N$. The classification in \cite{GIK} is based on the description of the automorphism group of a Kobayashi-hyperbolic Reinhardt domain obtained in \cite{Kr} (see also \cite{Sh}). The domains found in \cite{GIK} appear in part (i) of Theorem \ref{fullclasshom} as well as part (i) for $l=1$ and parts (ii)--(v) of Theorem \ref{main}. Further, in \cite{KV} simply connected {\it complete}\, Kobayashi-hyperbolic manifolds $X$ with $\dim G(X)=N^2$ were studied, where completeness means the completeness of the distance $K_X$. The main result of \cite{KV} states that every such manifold is biholomorphic to a Reinhardt domain and hence the classification in this case is a subset of that in \cite{GIK}.  Notice that the only manifolds on the list supplied by Theorem \ref{main} that are both simply connected and complete are those appearing in (ii). Together with the homogeneous domains in part (i) of Theorem \ref{fullclasshom} they form the partial classification obtained in \cite{KV}. At the same time, the statement of Theorem \ref{main} contains examples of Kobayashi-hyperbolic manifolds that cannot be realized as Reinhardt domains. For instance, it is not hard to show that the quotient of a spherical shell in (i) is not biholomorphic to any Reinhardt domain if $l>1$. 
\end{remark}

Next, the proof of Theorem \ref{main} can be modified to exclude real hypersurface orbits and yield a classification in the case $\dim G(X)=N^2-1$ for $N\ge 3$, with $X$ being non-homogeneous (see \cite{Isa2} and \cite{Isa6}, Theorem 4.1):

\begin{theorem}\label{mainn21} Let $N\ge 3$, $\dim G(X)=N^2-1$ and $X$ be non-homoge\-neous. Then $X$ is biholomorphic to $\BB^{N-1}\times S$, where $S$ is a hyperbolic Riemann surface with $\dim G(S)=0$.
\end{theorem}

It now remains to classify 2-dimensional Kobayashi-hyperbolic manifolds $X$ for which the group $G(X)$ is 3-dimensional. Determining such manifolds is the hardest part of our work, and we begin by giving a large number of examples (see \cite{Isa6}, Sect.~5.1 for details). Some of the examples are combined together, and some of them form parametric families. Note that the Kobayashi-hyperbolicity of the manifolds that appear below is not always obvious (see \cite{HI}). In what follows coordinates in $\CC^2$ are denoted by $(z,w)$.
\vspace{-0.1cm}\\

\noindent {\bf (i)} Here the group $G(X)$ consists of all maps of $\CC^2$ of the form
$$
\begin{array}{lll}
z & \mapsto & \lambda^{\alpha} z+i\beta,\\
w & \mapsto & \lambda w+i\gamma,\\
\end{array}
$$
where $\lambda>0$, $\beta,\gamma\in\RR$ and $\alpha\in\RR$ is a fixed number as specified below.
\vspace{-0.3cm}\\

{\bf (ia)} Fix $\alpha\in\RR$, $\alpha\ne 0,1$, and choose $s,t$ satisfying one of the following two sets of conditions:
\begin{itemize}

\item $0\le s<t\le\infty$, with either $s>0$ or $t<\infty$;

\item $-\infty\le s<0<t\le\infty$, where at least one of $s,t$ is finite, with $t\ne -s$ for $\alpha=1/2$.

\end{itemize}

Then let
$$
X=\left\{(z,w)\in\CC^2: s\left(\Re w\right)^{\alpha}<\Re z<t\left(\Re w\right)^{\alpha},\,\Re w>0\right\}.
$$
\vspace{-0.7cm}\\

{\bf (ib)} Fix $\alpha>0$, $\alpha\ne 1$, choose $-\infty<s<0<t<\infty$ and let
$$
\begin{array}{l}
X=\left\{(z,w)\in\CC^2: \Re z>s\left(\Re w\right)^{\alpha},\,\Re w>0\right\}\cup\\
\vspace{-0.3cm}\\
\hspace{0.95cm}\left\{(z,w)\in\CC^2:\Re z>t\left(-\Re w\right)^{\alpha},\,\Re w<0\right\}\cup\\
\vspace{-0.3cm}\\
\hspace{0.95cm}\left\{(z,w)\in\CC^2: \Re z>0,\, \Re w=0\right\}.
\end{array}
$$
\vspace{-0.5cm}\\

\noindent{\bf (ii)} Here the group $G(X)$ consists of all maps of $\CC^2$ of the form 
$$
\begin{array}{lll}
z & \mapsto & \lambda z+(\lambda\ln\lambda)w+i\beta,\\
w & \mapsto & \lambda w+i\gamma,\\
\end{array}
$$
where $\lambda>0$, $\beta,\gamma\in\RR$.
\vspace{-0.3cm}\\

{\bf (iia)} Choose $0\le s<t\le\infty$, with either $s>0$ or $t<\infty$, and let
$$
\begin{array}{l}
X=\left\{(z,w)\in\CC^2:\Re w\cdot\ln\left(s\Re w\right)<\Re z<\right.\\
\vspace{-0.3cm}\\
\hspace{6cm}\left.\Re w\cdot\ln\left(t\Re w\right),\,\Re w>0\right\}.
\end{array}
$$
\vspace{-0.7cm}\\

{\bf (iib)} Choose $-\infty<t<0<s<\infty$ and let
$$
\begin{array}{l}
X=\left\{(z,w)\in\CC^2:\Re z>\Re w\cdot\ln\left(s\Re w\right),\,\Re w>0\right\}\cup\\
\vspace{-0.3cm}\\
\hspace{0.95cm}\left\{(z,w)\in\CC^2:\Re z>\Re w\cdot\ln\left(t\Re w\right),\,\Re ,w<0\right\}\cup\\
\vspace{-0.3cm}\\
\hspace{0.95cm}\left\{(z,w)\in\CC^2: \Re z>0,\, \Re w=0\right\}.
\end{array}
$$
\vspace{-0.5cm}\\

\noindent {\bf (iii)} Here the group $G(X)$ consists of all maps of $\CC^2$ of the form
$$
\left(\hspace{-0.1cm}
\begin{array}{l}
z\\
w
\end{array}
\hspace{-0.1cm}\right)\mapsto A
\left(\hspace{-0.1cm}
\begin{array}{l}
z\\
w
\end{array}
\hspace{-0.1cm}\right)+i
\left(\hspace{-0.1cm}
\begin{array}{l}
\beta\\
\gamma
\end{array}
\hspace{-0.1cm}\right),
$$
where $A\in\SO_2(\RR)$, $\beta,\gamma\in\RR$.
\vspace{-0.3cm}\\

{\bf (iiia)} Choose $0\le s<t<\infty$ and let
$$
X=\left\{(z,w)\in\CC^2: s<\left(\Re z\right)^2+\left(\Re w\right)^2<t\right\}.
$$
\vspace{-0.7cm}\\

{\bf (iiib)} Choose $0<t<\infty$ and let
$$
X=\left\{(z,w)\in\CC^2: \left(\Re z\right)^2+\left(\Re w\right)^2<t\right\}.
$$
\vspace{-0.5cm}\\

\noindent {\bf (iv)} Here $X$ is any cover of any domain in {\bf (iiia)} (for details see \cite{Isa6}, pp.~64--66).
\vspace{-0.1cm}\\

\noindent{\bf (v)} Fix $\alpha>0$, choose $0<t<\infty$, $e^{-2\pi \alpha}t<s<t$ and let
$$
X=\left\{(z,w)\in\CC^2:se^{\alpha\phi}<r<te^{\alpha\phi}\right\},
$$  
where $(r,\phi)$ denote the polar coordinates in the $(\Re z,\Re w)$-plane with $\phi$ varying from $-\infty$ to $\infty$. The group $G(X)$ consists of all maps of $\CC^2$ of the form
$$
\left(\hspace{-0.1cm}
\begin{array}{l}
z\\
w
\end{array}
\hspace{-0.1cm}\right)\mapsto e^{\alpha\psi}
\left(\hspace{-0.1cm}
\begin{array}{rr}
\cos\psi & \sin\psi\\
-\sin\psi & \cos\psi
\end{array}
\hspace{-0.1cm}\right)
\left(\hspace{-0.1cm}
\begin{array}{l}
z\\
w
\end{array}
\hspace{-0.1cm}\right)+i
\left(\hspace{-0.1cm}
\begin{array}{l}
\beta\\
\gamma
\end{array}
\hspace{-0.1cm}\right),
$$
where $\psi,\beta,\gamma\in\RR$. 
\vspace{-0.1cm}\\

\noindent{\bf (vi)} Here the group $G(X)$ consists of all maps of $\CC\PP^2$ of the form
$$
(\zeta:z:w)\mapsto (\zeta:z:w) A,
$$
where $(\zeta:z:w)$ are homogeneous coordinates and $A\in\SO_3(\RR)$.
\vspace{-0.3cm}\\

{\bf (via)} Choose $1\le s<t<\infty$ and let
$$
\begin{array}{l}
X=\left\{(\zeta:z:w)\in\CC\PP^2: s |\zeta^2+z^2+w^2|<\right.\\
\vspace{-0.3cm}\\
\hspace{4cm}\left.|\zeta|^2+|z|^2+|w|^2< t |\zeta^2+z^2+w^2|\right\}.
\end{array}
$$
\vspace{-0.7cm}\\

{\bf (vib)} Choose $1<t<\infty$ and let
$$
X=\left\{(\zeta:z:w)\in\CC\PP^2: |\zeta|^2+|z|^2+|w|^2< t |\zeta^2+z^2+w^2|\right\}.
$$
\vspace{-0.7cm}\\

\noindent{\bf (vii)} Here $X$ is any cover of any domain in {\bf (via)} (for details see \cite{Isa6}, pp.~67--68). 
\vspace{-0.1cm}\\

\noindent{\bf (viii)} Choose $1<t<\infty$ and let
$$
X=\left\{(z_1,z_2,z_3)\in\CC^3:z_1^2+z_2^2+z_3^2=1,\, |z_1|^2+|z_2|^2+|z_3|^2<t \right\}.
$$
The group $G(X)$ consists of all maps of $\CC^3$ of the form
$$
\left(\hspace{-0.1cm}
\begin{array}{c}
z_1\\
z_2\\
z_3
\end{array}
\hspace{-0.1cm}\right)
\mapsto A
\left(\hspace{-0.1cm}
\begin{array}{c}
z_1\\
z_2\\
z_3
\end{array}
\hspace{-0.1cm}\right),
$$
where $A\in\SO_3(\RR)$.
\vspace{-0.1cm}\\

\noindent{\bf (ix)} Here the group $G(X)$ consists of all rational maps of $\CC^2$ of the form
$$
\left(\hspace{-0.1cm}
\begin{array}{c}
z\\
w
\end{array}
\hspace{-0.1cm}\right) \mapsto \displaystyle\frac{\left(\hspace{-0.1cm}\begin{array}{cc}
a_{11} & a_{12}\\
a_{21} & a_{22}
\end{array}
\hspace{-0.1cm}\right)\left(\hspace{-0.1cm}
\begin{array}{c}
z\\
w
\end{array}\hspace{-0.1cm}\right)+\left(\hspace{-0.1cm}
\begin{array}{c}
b_1\\
b_2
\end{array}
\hspace{-0.1cm}\right)}{c_1z+c_2w+d},
$$
where the matrix
\begin{equation}
Q:=\left(\hspace{-0.1cm}\begin{array}{ccc}
a_{11} & a_{12} & b_1\\
a_{21} & a_{22} & b_2\\
c_1& c_2 &d
\end{array}
\hspace{-0.1cm}\right)\label{matq}
\end{equation}
lies in $\SO_{2,1}^{\circ}$.
\vspace{-0.3cm}\\

{\bf (ixa)} Choose $-1\le s<t\le 1$ and let
$$
\begin{array}{l}
X=\left\{(z,w)\in\CC^2: s |z^2+w^2-1|<\right.\\
\vspace{-0.3cm}\\
\left.\hspace{4cm}|z|^2+|w|^2-1<t |z^2+w^2-1|\right\}.
\end{array}
$$
\vspace{-0.7cm}\\

{\bf (ixb)} Choose $-1<t<1$ and let
$$
X=\left\{(z,w)\in\CC^2: |z|^2+|w|^2-1< t |z^2+w^2-1|\right\}.
$$
\vspace{-0.5cm}\\

\noindent {\bf (x)} Here $X$ is any cover of any domain in {\bf (ixa)} (for details see \cite{Isa6}, pp.~70--74).
\vspace{-0.1cm}\\

\noindent{\bf (xi)} Here the group $G(X)$ consists of all rational maps of $\CC^2$ of the form
$$
\left(\hspace{-0.1cm}
\begin{array}{c}
z\\
w
\end{array}
\hspace{-0.1cm}\right) \mapsto \displaystyle\frac{\left(\hspace{-0.1cm}\begin{array}{cc}
a_{22} & b_2\\
c_2 & d
\end{array}
\hspace{-0.1cm}\right)\left(\hspace{-0.1cm}
\begin{array}{c}
z\\
w
\end{array}\hspace{-0.1cm}\right)+\left(\hspace{-0.1cm}
\begin{array}{c}
a_{21}\\
c_1
\end{array}
\hspace{-0.1cm}\right)}{a_{12}z+b_1w+a_{11}},
$$
where the matrix $Q$ defined as in (\ref{matq}) lies in $\SO_{2,1}^{\circ}$.
\vspace{-0.3cm}\\

{\bf (xia)} Choose $1\le s<t\le\infty$ and let
$$
\begin{array}{l}
X=\left\{(z,w)\in\CC^2:  s |1+z^2-w^2|<1+|z|^2-|w|^2<\right.\\
\vspace{-0.3cm}\\
\hspace{6cm}\left. t |1+z^2-w^2|,\,\Im\left(z(1+\overline{w})\right)>0\right\},
\end{array}
$$
where we assume that for $t=\infty$ the complex curve 
$$
\left\{(z,w)\in\CC^2: 1+z^2-w^2=0,\, \Im(z(1+\overline{w}))>0\right\}
$$
is not included in $X$.
\vspace{-0.3cm}\\

{\bf (xib)} Choose $1\le s<\infty$ and let
$$
\begin{array}{l}
X=\left\{(z,w)\in\CC^2:  1+|z|^2-|w|^2>s |1+z^2-w^2|,\right.\\
\vspace{-0.3cm}\\
\hspace{8cm}\left.\Im\left(z(1+\overline{w})\right)>0\right\}.
\end{array}
$$
\vspace{-0.5cm}\\

\noindent {\bf (xii)} Here $X$ is any cover of any domain in {\bf (xia)} (for details see \cite{Isa6}, pp.~70--74).
\vspace{-0.1cm}\\

\noindent{\bf (xiii)} Here $X$ is any manifold obtained by \lq\lq attaching\rq\rq, in a certain way, a Levi-flat and/or complex curve orbit to some of the manifolds introduced in {\bf (ix), (x), (xii)} (for details see \cite{Isa6}, pp.~74--79).
\vspace{-0.1cm}\\

\noindent{\bf (xiv)} $X=\BB^1\times S$, where $S$ is a hyperbolic Riemann surface with $\dim G(S)=0$ (cf.~Theorem \ref{mainn21}).
\vspace{-0.1cm}\\

We can now state our final result (see \cite{Isa3}; \cite{Isa6}, Chap.~5; \cite{Isa7}):

\begin{theorem}\label{class23} Let $N=2$ and $\dim G(X)=3$. Then $X$ is biholomorphic to one of the manifolds listed in {\bf (i)--(xiv)}.
\end{theorem}

Before closing the survey, we mention one significant application of results of this section. Namely, in \cite{V1}, \cite{V2} the classifications of 2-dimensional Kobayashi-hyperbolic manifolds $X$ with $\dim G(X)=3$ and $\dim G(X)=4$ were utilized to describe all domains $D$ in a complex surface $M$ such that: (i) $\overline{D}$ lies in a complete Kobayashi-hyperbolic domain in $M$ (e.g., $D$ is bounded if $M=\CC^2$), (ii) there exists an accumulation point $x\in\partial D$ of a $\Bihol(D)$-orbit, with $\partial D$ being real-analytic and of finite type near $x$. Problems related to determining domains with non-compact automorphism group are notoriously difficult (see, e.g., \cite{IKra1}), and the above result relies on our classification in an essential way. There are also other applications, but they deserve to be treated in a separate article, and we do not discuss them here.

\end{document}